\documentclass[12pt]{amsart}
\usepackage{a4wide,amsmath,enumerate,xcolor}
\usepackage{mathrsfs} % for \mathscr
\usepackage{esint} % for \fint
\allowdisplaybreaks
\usepackage[utf8]{inputenc}
\usepackage[english,vietnamese]{babel}

\let\pa\partial
\let\na\nabla
\let\eps\varepsilon

\newcommand{\R}{{\mathbb R}}
\newcommand{\T}{{\mathbb T}}
\newcommand{\diver}{\operatorname{div}}

\definecolor{dkgreen}{rgb}{0,0.45,0}

\newcommand{\dd}{{\mathrm{d}}}
\newcommand{\dx}{{\mathrm{d}x}}

\newtheorem{theorem}{Theorem}
\newtheorem{lemma}[theorem]{Lemma}
\newtheorem{proposition}[theorem]{Proposition}

\newtheorem{corollary}[theorem]{Corollary}
\newtheorem{definition}[theorem]{Definition}

\newtheorem{assumption}[theorem]{Assumption}

\begin{document}

\title[Stability for multi-species nonlocal cross-diffusion systems]{Weak--weak stability for multi-species \\
nonlocal cross-diffusion systems with \\ endpoint Lorentz kernels}

\author[A. J\"ungel]{Ansgar J\"ungel}
\address{Institute of Analysis and Scientific Computing, TU Wien, Wiedner Hauptstra\ss e 8--10, 1040 Wien, Austria}
\email{juengel@tuwien.ac.at} 

\author[\selectlanguage{vietnamese}T. T. Nguyễn]{ \selectlanguage{vietnamese}Tuấn Tùng Nguyễn}
\address{Institute of Analysis and Scientific Computing, TU Wien, Wiedner Hauptstra\ss e 8--10, 1040 Wien, Austria}
\email{tung.nguyen@tuwien.ac.at}

\author[Z. Sun]{Zhiwei Sun}
\address{Department of Mathematics, Hong Kong University of Science and Technology, Clear Water Bay, Kowloon, Hong Kong}
\email{zwsun@ust.hk} 

\date{\today}
\selectlanguage{english}

\thanks{The first two authors acknowledge partial support from the Austrian Science Fund (FWF), grant 10.55776/PAT2687825, and from the Austrian Federal Ministry for Women, Science and Research and implemented by \"OAD, project MULT09/2025. This work has received funding from the European Research Council (ERC) under the European Union's Horizon 2020 research and innovation programme, ERC Advanced Grant NEUROMORPH, no.~101018153. For open-access purposes, the authors have applied a CC BY public copyright license to any author-accepted manuscript version arising from this submission.} 

\begin{abstract}
The quantitative weak--weak stability for multi-species nonlocal cross-diffu\-sion systems on the torus is established. The interactions of the species are described by a symmetric, positive semidefinite matrix-valued kernel whose gradient belongs to the endpoint Lorentz space $L^{d,\infty}$, where $d\ge 2$ is the space dimension. The central estimate, formulated in terms of a modulated Hellinger distance--interaction energy, yields uniqueness of solutions with finite Boltzmann entropy and Fisher information. The proof relies on a two-level endpoint mechanism: The Boltzmann entropy controls the interaction energy, while the Fisher information controls the associated force field. 
%The results extend to bounded interaction kernels for $d\ge1$, without requiring regularity of their derivatives.
Further results are the immediate $L^2$ regularization and the exponential decay to the homogeneous equilibrium. As a particular application, uniqueness and stability for the two-dimensional Poisson--Nernst--Planck system with periodic Coulomb kernels are obtained.
\end{abstract}

% \paragraph{Keywords:}  
 
% \paragraph{AMS classification:}  

\maketitle

%%%%%%%%%%%%%%%%%%%%%%%%%%%%%%%%%%%%%%%%%%%%%%%%%%%%%%%%%%%%%%%%%%%

\section{Introduction}

Nonlocal cross-diffusion systems model the evolution of densities that diffuse while being transported by nonlocal interaction fields. Such systems arise in interacting particle models \cite{CDJ19}, aggregation--diffusion systems \cite{CSS26}, multi-species population dynamics \cite{JPZ22}, and mean-field network models \cite{BFFT12}. While smooth interaction kernels have been studied intensively in the literature (see, e.g., \cite{GHLP22}), fewer works address irregular kernels. In particular, uniqueness within weak solution classes has remained open. In this paper, we resolve this problem in the repulsive regime and, moreover, establish a quantitative weak--weak stability result.

We consider the multi-species nonlocal cross-diffusion equations
\begin{align}\label{1.eq}
  \pa_t u_i = \Delta u_i + \diver(u_i\na p_i[u]), \quad  
  p_i[u] = \sum_{j=1}^n K_{ij}*u_j, \quad i=1,\ldots,n,
\end{align}
on the flat torus $\T^d$ ($d\ge 1$), with the initial conditions
\begin{align}\label{1.ic}
  u_i(0) = u_{i,0}\ge 0\quad \mbox{in }\T^d,\ i=1,\ldots,n.
\end{align}
Here, $u=(u_1,\ldots,u_n)$ denotes the vector of the species densities, and $K=(K_{ij})_{i,j=1}^n$ is a matrix of periodic interaction kernels. The model describes the nonlocal interaction of the species, as organisms (or agents) often respond to stimuli originating beyond their immediate local environment. The ability of organisms to detect a given stimulus and respond accordingly is encoded in the interaction kernels. Equations \eqref{1.eq} can be derived from interacting particle systems in the many-particle limit; see, e.g., \cite{CDJ19}. When the kernels converge to the Dirac delta distribution, up to some factors, we recover the Busenberg--Travis model  with additional diffusion $\Delta u_i$ \cite{BuTr83}. In this paper, we consider kernels with possible singularities and aim to determine the critical singularity scale at which the specific weak solution structure is sufficient to close the stability argument.

Before making this statement more explicit, let us briefly review the literature. For smooth interaction kernels, local well-posedness results for the nonlocal cross-diffusion system \eqref{1.eq} were shown in \cite{GHLP22}. The existence of global weak solutions to \eqref{1.eq} is proved in \cite{JPZ22} for positive semidefinite kernels satisfying a detailed-balance condition. In particular, weak--strong uniqueness is proved for $L^2(0,T;H^1(\T^d))$ solutions when $\na K_{ij}\in L^{d+2}$. Global results were obtained in \cite{DiFa13} for kernels with a singularity at the origin. More recently, \cite{CSS26} developed a well-posedness theory for the aggregation dynamics with irregular bounded kernels such that $K_{ij}\in \mbox{BV}\cap L^\infty$, i.e., whose distributional gradients are finite measures, and a uniqueness result was established for strong solutions. The authors of \cite{CGH25} prove the global existence of weak solutions for small initial data and interaction kernels potentially as singular as the Newtonian potential. Related works investigated smooth solutions under smallness assumptions \cite{CoRa26} and heterogeneous systems \cite{MMR26}. The nonlocal-to-local limit was shown in \cite{DHPP24}, even without additional diffusion. These works address important low-regularity and singular-kernel regimes, but the question of uniqueness within weak solution classes remained open. The present paper addresses this issue.

In this work, a key structural assumption is that the interaction kernel 
$K=(K_{ij})_{i,j=1}^n$ is positive semidefinite, in the sense that the interaction energy
\begin{equation}\label{eq:def-interaction-energy}
\mathcal P[u]
:=
\frac12
\sum_{i,j=1}^n
\int_{\T^d}
u_i\,K_{ij}*u_j\,\dd x
\end{equation}
is nonnegative for any admissible functions (see Assumption~\ref{ass.kern}). In this case, equations \eqref{1.eq} admit an entropy structure and formally, the entropy dissipation identity
\begin{align}\label{1.ent}
  \frac{\dd}{\dd t}\sum_{i=1}^n\int_{\T^d}u_i(\log u_i-1)\dx 
  + 4\sum_{i=1}^n\int_{\T^d}|\na\sqrt{u_i}|^2\dx 
  + \mathcal D[u] = 0
\end{align}
holds along smooth solutions to \eqref{1.eq}, where the dissipation $\mathcal D[u]$ is given by
\begin{align}\label{1.D}
  \mathcal D[u] = \sum_{i=1}^n\int_{\T^d}
  \na u_i\cdot\na p_i[u]\dx.
\end{align}
Moreover, $\mathcal D[u]\ge0$ by the positive semidefiniteness of $K$. This suggests seeking solutions  
satisfying 
\begin{align}\label{1.weak}
  u_i\ge 0, \quad \int_{\T^d}u_i\log u_i\dx < \infty, \quad
  \|\na\sqrt{u_i}\|_{L^2((0,T)\times\T^d)}<\infty
\end{align}
for all $i=1,\ldots,n$, where the last two conditions correspond to finite Boltzmann entropy and Fisher information. We refer to this regularity class as the entropy--Fisher class. We show that this class is sufficient for weak--weak uniqueness and stability, without any a priori $L^\infty(\T^d)$ bounds or $L^2(\T^d)$-based stability of the solution.

For $d\ge2$, let $K=(K_{ij})$ be a symmetric, periodic, positive semidefinite matrix of kernels satisfying
\begin{align}\label{1.K}
  \na K_{ij}\in L^{d,\infty}(\T^d;\R^d)\quad\mbox{for }i,j=1,\ldots,n.
\end{align}
Here, the Lorentz space $L^{d,\infty}(\T^d)$ is a generalization of the standard $L^{p}(\T^d)$ spaces; see Section \ref{sec.lor}. The main result of this paper is that
any two entropy--Fisher weak solutions $u$ and $v$ to \eqref{1.eq} (i.e.\ solutions satisfying \eqref{1.weak}) with initial data $u_0$ and $v_0$, respectively, satisfy for some $\lambda>0$ the quantitative stability estimate
\begin{align}\label{1.stab}
  \mathcal H[u|v](t) + \lambda\mathcal P[u-v](t)
  \le C_T\big(\mathcal H[u_0|v_0] + \lambda\mathcal P[u_0-v_0]\big),
\end{align}
where the (squared) Hellinger distance $\mathcal{H}$ is defined by
\begin{align}\label{1.HP}
  \mathcal H[u|v] = \sum_{i=1}^n\int_{\T^d}
  \big(\sqrt{u_i}-\sqrt{v_i}\big)^2\dx,
\end{align}
and the interaction energy $\mathcal P$ is given by \eqref{eq:def-interaction-energy}. If the initial data coincide, $u_0=v_0$, inequality \eqref{1.stab} implies that $u(t)=v(t)$ in $\T^d$ for $t>0$, and weak--weak uniqueness of solutions follows. To the best of our knowledge, the coupling of the Hellinger distance and interaction energy has not yet been explored in the literature.

We claim that condition \eqref{1.K} is the critical endpoint scale at which the entropy--Fisher regularity \eqref{1.weak}, generated by the diffusion, matches the singularity of the nonlocal interaction. Indeed, at the level of the force field $\na p_i[u]$, it follows from the endpoint Sobolev estimate and the O'Neil--Young convolution inequality that
\begin{align}\label{1.nap}
  \|\na p_i[u]\|_{L^2(0,T;L^\infty(\T^d))}
  \le C\sum_{j=1}^n\sqrt{m_j}\|\na K_{ij}\|_{L^{d,\infty}(\T^d)}
  \|\na\sqrt{u_j}\|_{L^2((0,T)\times\T^d)},
\end{align}
where $C>0$ and $m_j$ denotes the mass of $u_j$; see Lemma \ref{lem.ec2} below. Thus, the Fisher information $\|\na\sqrt{u_i}\|_{L^2((0,T)\times\T^d)}$ controls the force field. At the level of the interaction energy $\mathcal P$, condition \eqref{1.K} implies that $K_{ij}\in\mbox{BMO}(\T^d)$, where BMO is the space of functions of bounded mean oscillation, introduced in Section \ref{sec.lor}. Together with the $L\log L$--BMO pairing, this gives
\begin{align}\label{1.P}
  \mathcal P[u] \le C\sum_{i,j=1}^n m_j\|K_{ij}\|_{\text{BMO}(\T^d)}
  \|u_i\|_{L\log L(\T^d)},
\end{align}
where $C>0$; see Lemma \ref{lem.ec1}. The right-hand side is bounded in view of \eqref{1.weak}. This shows that the entropy in $L\log L$ controls the interaction energy, while the Fisher information controls the force field. This two-level control is key to closing
the weak--weak stability estimate.

Our result strengthens the existing literature in two directions. First, it upgrades the $L^2(\T^d)$-based weak--strong uniqueness of \cite{JPZ22} to weak--weak uniqueness in the entropy--Fisher class. Second, it lowers the kernel regularity to the critical endpoint regime \eqref{1.K}, which includes important singular examples, most notably the two-dimensional Coulomb kernel (see Section \ref{sec.PNP}). The same results also hold for the alternative assumption of bounded kernels $K_{ij}\in L^\infty(\T^d)$ for $d\ge1$, without regularity conditions on their derivatives. In this case, Young's
convolution inequality provides simpler
estimates for the force field and interaction energy, replacing the endpoint Lorentz estimates and the $L\log L$--BMO pairing.

A further consequence of our approach is that the unique entropy--Fisher weak solution becomes an $L^2(\T^d)$ function at a.e.\ positive time, thus improving the regularity from $L\log L(\T^d)$ to $L^2(\T^d)$. Moreover, we conclude the exponential decay
\begin{align*}
  \sum_{i=1}^n m_i^{-1}\|u_i(t)-m_i\|_{L^2(\T^d)}^2
  \le C\big(e^{2\kappa t/d}-1\big)^{-d/2} \quad\mbox{for a.e. }t>0,
\end{align*}
where $\kappa=4/C_{LS}$ and $C_{LS}>0$ is the constant of the logarithmic Sobolev inequality; see Lemma \ref{lem.lsi}. This profile combines heat-type $L^1(\T^d)$-to-$L^2(\T^d)$
regularization at short times with exponential relaxation at large times.

Before stating our main results, we outline the key ideas behind the stability estimate \eqref{1.stab} in the scalar case. The same principle extends to the multi-species situation. Let $u$ and $v$ be two smooth positive solutions to \eqref{1.eq} with $n=1$. %Let $L=\log(u/v)$ and $W=\na(p[u]-p[v])$. 
Then, by Young's inequality,
\begin{align}\label{1.estH}
  \frac{\dd}{\dd t}&\mathcal H[u|v]
  + \frac12\int_{\T^d} \sqrt{uv}\bigg|\na\log\frac{u}{v}\bigg|^2\dx
  = -\frac12\int_{\T^d}\sqrt{uv}\na(p[u]-p[v])
  \cdot\na\log\frac{u}{v}\dx \\
  &\le \frac14\int_{\T^d}\sqrt{uv}\bigg|\na\log\frac{u}{v}\bigg|^2\dx
  + \frac14\int_{\T^d}\sqrt{uv}|\na(p[u]-p[v])|^2\dx \nonumber \\
  &\le \frac14\int_{\T^d}\sqrt{uv}\bigg|\na\log\frac{u}{v}\bigg|^2\dx
  + \frac18\int_{\T^d}(u+v)|\na(p[u]-p[v])|^2\dx, \nonumber
\end{align}
where $p[u]=K*u$ and $\mathcal H[u|v]$ is given in \eqref{1.HP}. To close this estimate, we compute
\begin{align*}
  \frac{\dd}{\dd t}&\mathcal P[u-v] + \mathcal D[u-v]
  + \frac12\int_{\T^d}(u+v)|\na(p[u]-p[v])|^2\dx \\
  &= -\frac12\int_{\T^d}(u-v)\na(p[u]+p[v])
  \cdot\na(p[u]-p[v])\dx \nonumber \\
  &\le \frac14\int_{\T^d}(u+v)|\na(p[u]-p[v])|^2\dx
  + \frac14\|\na(p[u]+p[v])\|_{L^\infty(\T^d)}^2
  \int_{\T^d}\frac{(u-v)^2}{u+v}\dx, \nonumber 
\end{align*}
recalling definitions \eqref{1.D} and \eqref{eq:def-interaction-energy} of $\mathcal D$ and $\mathcal P$. The first term on the right-hand side is absorbed by the left-hand side, while we use $(u-v)^2/(u+v)\le 2(\sqrt{u}-\sqrt{v})^2$ for the second term. This shows that
\begin{align}\label{1.estP}
  \frac{\dd}{\dd t}\mathcal P[u-v] &+ \mathcal D[u-v]
  + \frac14\int_{\T^d}(u+v)|\na(p[u]-p[v])|^2\dx \\
  &\le \frac12\|\na(p[u]+p[v])\|_{L^\infty(\T^d)}^2\mathcal H[u|v]. \nonumber 
\end{align}
We add \eqref{1.estH} and $\lambda$ times \eqref{1.estP} for some $\lambda>0$ and neglect some nonnegative contributions to find that 
\begin{align*}
  \frac{\dd}{\dd t}(\mathcal H[u|v]&+\lambda\mathcal P[u-v])
  + \bigg(\frac{\lambda}{4}-\frac18\bigg)
  \int_{\T^d}(u+v)|\na(p[u]-p[v])|^2\dx \\
  &\le \frac{\lambda}{2}\|\na(p[u]+p[v])\|_{L^\infty(\T^d)}^2
  \mathcal H[u|v].
\end{align*}
Adding $(\lambda^2/2)\|\na(p[u]+p[v])\|_{L^\infty(\T^d)}^2\mathcal P[u-v]\ge 0$ to the right-hand side and choosing $\lambda\ge 1/2$, we can apply Gronwall's lemma to derive the stability estimate \eqref{1.stab}, provided that $\na p[u]$, $\na p[v]\in L^2(0,T;L^\infty(\T^d;\R^d))$. The kernel assumption \eqref{1.K} gives exactly this regularity in the entropy--Fisher class thanks to \eqref{1.nap}. 
The argument above is only formal for
entropy--Fisher weak solutions. Indeed,
the densities may vanish, and their limited
time regularity does not allow us to justify
the evolution of the modulated energies
directly. To overcome these difficulties,
we compare an arbitrary weak solution with
a solution of a kernel-regularized system.
The rigorous argument combines a renormalized
formulation for the Hellinger distance with
a regularized interaction energy identity,
while the endpoint estimates provide the
control needed to pass to the limit.

The paper is organized as follows. We state our main assumptions and results in Section \ref{sec.main}. We recall basic functional inequalities and the definitions of Lorentz and BMO spaces, together with the proof of the endpoint estimates \eqref{1.nap} and \eqref{1.P} in Section \ref{sec.prep}. In Section \ref{sec.ex}, we prove the existence of an entropy--Fisher weak solution to \eqref{1.eq}--\eqref{1.ic} by passing to the limit in approximate solutions to \eqref{1.eq} with regularized kernels. The Hellinger estimate related to \eqref{1.estH} is proved in Section \ref{sec.hell}, while the relative interaction energy estimate \eqref{1.estP} and the weak--weak stability estimate \eqref{1.stab} are shown in Section \ref{sec.ener}. Section \ref{sec.regul} is devoted to the proof of the $L^2(\T^d)$ regularization and the convergence to equilibrium. Finally, Section~\ref{sec:bounded-kernels} extends the analysis to bounded interaction kernels and provides examples of admissible kernels.

%%%%%%%%%%%%%%%%%%%%%%%%%%%%%%%%%%%%%%%%%%%%%%%%%%%%%%%%%%%%%%%%%%%

\section{Main results}\label{sec.main}

We prove three main results in this paper: the global existence of entropy--Fisher weak solutions to \eqref{1.eq}--\eqref{1.ic}, their  weak--weak modulated stability, and the positive-time regularization and convergence to equilibrium. The results are applied to a repulsive Poisson--Nernst--Planck system in two space dimensions. First, we specify our assumptions. 

\begin{assumption}[Matrix kernel]\label{ass.kern}\rm 
The flat torus is defined by $\T^d=\R^d/\mathbb Z^d$ with $|\T^d|=1$. The kernels $K_{ij}\in L^1(\T^d)$ are periodic and satisfy $\na K_{ij}\in L^{d,\infty}(\T^d;\R^d)$, the symmetry assumption
\begin{align*}
  K_{ij}(x) = K_{ji}(-x)\quad\mbox{for a.e. }x\in\T^d,\
  i,j=1,\ldots,n,
\end{align*}
the zero-mean condition $\int_{\T^d}K_{ij}(x)\dx = 0$, and the matrix kernel $K=(K_{ij})$ is positive semidefinite in the sense
\begin{align*}
  \sum_{i,j=1}^n\int_{\T^d}u_i(K_{ij}*u_j)\dx \ge 0
  \quad\mbox{for all }u=(u_1,\ldots,u_n)\in L^2(\T^d;\R^n).
\end{align*}
\end{assumption}

\begin{assumption}[Initial data]\label{ass.init}\rm 
The initial datum $u_0=(u_{1,0},\ldots,u_{n,0})$ is componentwise nonnegative, and satisfies $u_{i,0}\log u_{i,0}\in L^1(\T^d)$, $m_i:=\int_{\T^d}u_{i,0}\dx>0$ for $i=1,\ldots,n$. 
\end{assumption}

The positive semidefiniteness in Assumption~\ref{ass.kern}
also implies the nonnegativity of the dissipation \eqref{1.D}, since
\begin{align*}
\mathcal D[u]
=2\sum_{k=1}^d\mathcal P[\partial_{x_k}u]\ge0
\quad\mbox{for all }u=(u_1,\ldots,u_n)\in H^1(\T^d;\R^n)
\end{align*}

\subsection{Entropy--Fisher weak solutions and global existence}

We introduce the following notion of weak solution. Let $T>0$.

\begin{definition}[Entropy--Fisher weak solution]\label{def.ex}
A componentwise nonnegative function $u=(u_1,\ldots,u_n)$ is called an {\em entropy--Fisher weak solution} to \eqref{1.eq}--\eqref{1.ic} on $(0,T)\times\T^d$ if the mass is conserved, $\|u_i(t)\|_{L^1(\T^d)}=m_i$ for a.e.\ $t\in(0,T)$, the regularity
\begin{align}\label{1.bd}
  & u_i(t)\log u_i(t)\in L^1(\T^d)\text{ for a.e.}\ t\in(0,T), \quad
  \na\sqrt{u_i}\in L^2(0,T;L^2(\T^d;\R^d)), \\
  & \na u_i+u_i\na p_i[u]\in L^2(0,T;L^1(\T^d;\R^d)), \quad
  \pa_t u_i \in L^2(0,T;W^{-1,1}(\T^d)) \nonumber 
\end{align}
holds, where $W^{-1,1}(\T^d)$ is the dual space of $W^{1,\infty}(\T^d)$, and for any $\varphi\in W^{1,\infty}((0,T)\times\T^d)$ with $\varphi(T)=0$ and $i=1,\ldots,n$, we have the weak formulation
\begin{align*}
  -\int_0^T\int_{\T^d}u_i\pa_t\varphi \dx\dd t
  - \int_{\T^d}u_{i,0}\varphi(0)\dx
  = -\int_0^T\int_{\T^d}(\na u_i+u_i\na p_i[u])\cdot\na\varphi \dx\dd t.
\end{align*}
\end{definition}

\begin{proposition}[Existence of entropy--Fisher weak solutions]
\label{prop.ex} 
Let $d\ge2$. Let Assumptions \ref{ass.kern} and \ref{ass.init} hold. Then there exists a global entropy--Fisher weak solution to \eqref{1.eq}--\eqref{1.ic} in the sense of Definition \ref{def.ex}.
\end{proposition}

The proof is based on a regularization strategy. With $K^\delta$ being a regularized kernel, we show first the existence of a solution $u^\delta$ to \eqref{1.eq} with $K$ replaced by $K^\delta$. Uniform estimates are derived from the entropy identity \eqref{1.ent}, and compactness is obtained from the Aubin--Lions lemma.

%%%%%%%%%%%%%%%%%%%%%

\subsection{Weak--weak modulated stability and uniqueness}
%\label{sec.stab}

This section contains our main result. Recall definition \eqref{1.D} of the dissipation $\mathcal D$ and definitions \eqref{1.HP} and \eqref{eq:def-interaction-energy} of the Hellinger distance $\mathcal H$ and the interaction energy $\mathcal P$, respectively. 

\begin{theorem}[Weak--weak modulated stability]\label{thm.stab}
Let $d\ge2$. Let Assumption \ref{ass.kern} hold and let $u$ and $v$ be two entropy--Fisher weak solutions to \eqref{1.eq}--\eqref{1.ic} with initial data $u_0$ and $v_0$, respectively, satisfying Assumption \ref{ass.init}. Then there exist constants $\lambda>0$, $c>0$, and $C_T>0$, only depending on $T$ as well as $d$, $n$, the mass vector, the endpoint kernel norms, and the entropy--Fisher bounds \eqref{1.bd} of $u$ and $v$, such that for a.e.\ $t\in(0,T)$,
\begin{align*}
  \mathcal M_\lambda[u(t)|v(t)] + c\int_0^t\mathcal D[u(s)-v(s)]\dd s
  \le C_T\mathcal M_\lambda[u_0|v_0],
\end{align*} 
where $\mathcal M_\lambda[u|v] = \mathcal H[u|v]+\lambda\mathcal P[u-v]$ is the modulated Hellinger distance--interaction energy. Moreover, if $u_0=v_0$, then $u(t)=v(t)$ in $\T^d$ for a.e.\ $t\in(0,T)$, i.e., entropy--Fisher weak solutions are unique. 
\end{theorem}

As explained in the introduction, the proof relies on the derivation of inequalities \eqref{1.estH} and \eqref{1.estP}, which is carried out at the level of the approximate solutions $v^\delta$. First, we derive a renormalized formulation of \eqref{1.eq} and regularize the square root $\sqrt{u_i}$ by $\sqrt{u_i+\eps}$ for $\eps>0$. This allows us to derive the Hellinger estimate \eqref{1.estH} for $u$ and $v^\delta$. Second, we estimate the interaction energy by estimating the relative field
\begin{align*}
  \na p_i[u]-\na p_i^\delta[v^\delta]
  = \na(p_i^\delta[u]-p_i^\delta[v^\delta]) + E_i^\delta[u], \quad
  E_i^\delta[u] := \na(p_i[u]-p_i^\delta[u]),
\end{align*}
where $p_i^\delta[u]$ is defined as in \eqref{1.eq} with $K_{ij}$ replaced by $K_{ij}^\delta$. This leads to an approximate version of \eqref{1.estP}, with $E_i^\delta[u]\to 0$ strongly in $L^2(0,T;L^\infty(\T^d))$ as $\delta\to 0$. Combining \eqref{1.estH} and \eqref{1.estP} then proves the weak--weak stability estimate. 

%%%%%%%%%%%%%%%%%%%%

\subsection{Positive-time regularization and convergence to equilibrium}
%\label{sec.time}

We show that the unique entropy--Fisher weak solution is dynamically self-improving: It instantaneously enters the stronger $L^2(\T^d)$ regime and converges exponentially fast to equilibrium.

\begin{theorem}[Positive-time regularization and decay]\label{thm.time}
Let $d\ge2$. Let Assumptions \ref{ass.kern} and \ref{ass.init} hold and let $u$ be the unique entropy--Fisher weak solution to \eqref{1.eq}--\eqref{1.ic}. Then, for any $\tau>0$, $u\in L^\infty(\tau,\infty;L^2(\T^d))$ and there exists $C>0$ independent of $t$, such that
\begin{align*}
  \sum_{i=1}^n\frac{1}{m_i}\|u_i(t)-m_i\|_{L^2(\T^d)}^2
  \le C(e^{2\kappa t/d}-1)^{-d/2}\quad\mbox{for a.e. }t>0,
\end{align*}
where $m_i=\int_{\T^d}u_{i,0}\dx$ and $\kappa>0$ only depends on $d$.
\end{theorem}

Since $\left(e^{2\kappa t/d}-1\right)^{-d/2}
\sim C t^{-d/2}$ as $t\to0$, the short-time smoothing rate in Theorem~\ref{thm.time} is the same as that of the linear heat equations.
The exponential decay in the $L^1(\T^d)$ norm follows from the relative entropy, logarithmic Sobolev, and Csisz\'ar--Kullback--Pinsker inequalities (see Section \ref{sec.prep}). The main difficulty lies in establishing a quantitative
$L^2(\T^d)$ regularization and decay estimate uniformly,
since the initial data need not belong to $L^2(\T^d)$. To this end, we estimate the weighted $L^2(\T^d)$ norm of $\na w_i^\delta$, where $w_i^\delta=u_i^\delta-m_i$ and $u^\delta$ is a regularized solution, using the zero-mean Nash inequality. This yields a differential inequality for a suitable weighted $L^2(\T^d)$ norm, which can be solved explicitly and leads to the decay rate stated in Theorem \ref{thm.time}. 

%%%%%%%%%%%%%%%%%%%%

\subsection{Two-dimensional Poisson--Nernst--Planck systems}\label{sec.PNP}

As a concrete example of the abstract kernel framework, we consider the Poisson--Nernst--Planck system on the two-dimensional torus,
\begin{align}\label{1.PNP}
  \pa_t u_i = \Delta u_i + \diver(z_iu_i\na\phi[u]), \quad
  -\Delta\phi[u] = \rho[u] - \int_{\T^2}\rho[u]\dx\quad\mbox{in }\T^2,
  \ t>0,
\end{align}
where $\rho[u]=\sum_{j=1}^n z_ju_j$ is the total charge density. Here, $u_i\ge 0$ are particle densities, $z_i\in\R$ are their charges, and $\phi[u]$ is an electric potential with zero mean. This system is a special case of \eqref{1.eq} with $p_i[u]=z_i\phi[u]$ and $K_{ij}=z_iz_jG$, where $G$ is the zero-mean Green function of the periodic Laplacian, defined by
\begin{align*}
  -\Delta G = \delta_0-1\quad\mbox{in }\T^2, \quad
  \int_{\T^2}G(x)\dx = 0.
\end{align*}
The periodic Green function satisfies $\na G\in L^{2,\infty}(\T^2;\R^2)$ and consequently $\na K_{ij}\in L^{2,\infty}(\T^2)$. Note that the periodic Coulomb kernel is unbounded, and hence the
theory of \cite{CSS26}, which requires $K_{ij}\in \mbox{BV}\cap L^\infty$, does not apply. Moreover, $K$ is symmetric and positive semidefinite, since for any $u=(u_1,\ldots,u_n)\in L^2(\T^2;\R^n)$,
\begin{align*}
  \sum_{i,j=1}^n\int_{\T^2}u_i(K_{ij}*u_j)\dx
  = \int_{\T^2}\rho[u](G*\rho[u])\dx
  = \int_{\T^2}\rho[u]\phi[u]\dx 
  = \int_{\T^2}|\na\phi[u]|^2\dx\ge 0.
\end{align*}
Thus, Assumption \ref{ass.kern} is satisfied, and we conclude the following result.

\begin{corollary}[Weak--weak stability and uniqueness of the Poisson--Nernst--Planck model]
Let $d=2$. Entropy--Fisher weak solutions to \eqref{1.PNP} with initial condition \eqref{1.ic} and same mass are unique and satisfy
\begin{align*}
  \mathcal H[u|v](t) + \lambda\mathcal P_{\rm el}[u-v](t)
  + c\int_0^t\int_{\T^2}\rho[u-v]^2\dx\dd s
  \le C_T\big(\mathcal H[u_0|v_0] + \lambda
  \mathcal P_{\rm el}[u_0-v_0]\big)
\end{align*}
for a.e.\ $t\in(0,T)$, where
\begin{align*}
  \mathcal P_{\rm el}[u-v] = \frac12\int_{\T^2}|\na(\phi[u]-\phi[v])|^2
  \dx
\end{align*}
is the electric energy. Moreover, the unique entropy--Fisher weak solution satisfies the positive-time $L^2(\T^2)$ regularization and decay estimates of Theorem \ref{thm.time}.
\end{corollary}

\subsection{Bounded interaction kernels}

The results above are established for $d\ge2$,
which is needed by the endpoint Lorentz estimates. We now show that the same conclusions hold in any space dimension $d\ge1$
under an alternative kernel assumption.
Instead of the
endpoint Lorentz condition on $\nabla K_{ij}$, we assume that the kernels
$K_{ij}$ are bounded.

\begin{theorem}[Bounded interaction kernels]
\label{thm:bounded-kernels}
Let $d\ge1$. Suppose that $K_{ij}\in L^\infty(\T^d)$ for
$i,j=1,\ldots,n$, and that the symmetry, zero-mean, and positive
semidefiniteness conditions in Assumption~\ref{ass.kern} hold.
Let the initial data satisfy Assumption~\ref{ass.init}.
Then the conclusions of
Proposition~\ref{prop.ex},
Theorem~\ref{thm.stab}, and
Theorem~\ref{thm.time}
remain valid.
\end{theorem}

The bounded-kernel assumption allows for kernels without regularity assumptions on their derivatives. The proof relies on the same framework
as above, with simpler estimates for the drift field and interaction energy based on Young's convolution inequality.
In particular, the interaction energy is controlled by the $L^1$ norms of the densities, without using
the $L\log L$--BMO pairing. The remaining arguments carry over as explained in
Section~\ref{sec:bounded-kernels},
where we also present
some examples of admissible bounded kernels.

%%%%%%%%%%%%%%%%%%%%%%%%%%%%%%%%%%%%%%%%%%%%%%%%%%%%%%%%%%%%%%

\section{Preparations}\label{sec.prep}

We recall some functional inequalities, the definition of Lorentz and BMO spaces, and prove some endpoint estimates.

\subsection{Basic functional inequalities}

We recall some functional inequalities. The first two inequalities are stated in \cite[Sec.~2.2, Appendix A9]{Jue16}, while the third one is due to Nash \cite{Nas58}.

\begin{lemma}[Logarithmic Sobolev inequality]\label{lem.lsi}
There exists a constant $C_{LS}>0$, only depending on $d$, such that for all $f\in H^1(\T^d)$ with $f\not\equiv 0$,
\begin{align*}
  \int_{\T^d}f^2\log\frac{f^2}{\|f\|_{L^2(\T^d)}^2}\dx
  \le C_{LS}\int_{\T^d}|\na f|^2 \dx. 
\end{align*}
\end{lemma}

\begin{lemma}[Csisz\'ar--Kullback--Pinsker inequality]\label{lem.ckpi}
It holds for all nonnegative functions $u\in L^1(\T^d)$ with $m=\int_{\T^d}u\dx>0$ that
\begin{align*}
  \|u-m\|_{L^1(\T^d)}^2 \le 2m\int_{\T^d}u\log\frac{u}{m}\dx.
\end{align*}
\end{lemma}

\begin{lemma}[Zero-mean Nash inequality]\label{lem.nash}
There exists a constant $C_N>0$, depending only on $d$, such that for all $f\in H^1(\T^d)$ with $\int_{\T^d}f\dx=0$,
\begin{align*}
  \|f\|_{L^2(\T^d)}^{2+4/d}\le C_N\|\na f\|_{L^2(\T^d)}^2
  \|f\|_{L^1(\T^d)}^{4/d}.
\end{align*}
\end{lemma}

%%%%%%%%%%%%%%%%

\subsection{Lorentz and BMO spaces}\label{sec.lor}

For any function $f\in L^1(\T^d)$, we denote its mean value by $\langle f\rangle = |\T^d|^{-1}\int_{\T^d}f\dx$. Let $1\le p<\infty$ and $1\le q\le\infty$. The Lorentz space $L^{p,q}(\T^d)$ is the space of all measurable functions $\T^d\to\R$ with finite quasinorm
\begin{align*}
  \|f\|_{L^{p,q}(\T^d)} &= \bigg(\int_0^\infty(t^{1/p}
  f^*(t))^q\frac{\dd t}{t}\bigg)^{1/q} \quad\mbox{if }q<\infty, \\
  \|f\|_{L^{p,\infty}(\T^d)} &= \sup_{t>0}t^{1/p}f^*(t)
  \quad\mbox{if }q=\infty,
\end{align*}
where $f^*$ is the nonincreasing rearrangement of $|f|$,
\begin{align*}
  f^*(t) = \inf\big\{\lambda>0:|\{x\in\T^d:|f(x)|>\lambda\}|\le t
  \big\}, \quad t>0.
\end{align*}
The Lorentz spaces are complete quasinormed spaces, and $L^{p,\infty}(\T^d)$ is the weak Lebesgue space. It holds that $L^{p,p}(\T^d)=L^p(\T^d)$ and $L^{p,q}(\T^d)\subset L^{p,r}(\T^d)$ for $q\le r$ \cite[Sec.~1.4.2]{Gra14}. This implies that $L^{p,q}(\T^d)\subset L^{p,p}(\T^d)=L^p(\T^d)$ for $q\le p$. When $q>p$, we have only $L^{p,q}(\T^d)\subset L^s(\T^d)$ for all $s<p$. In particular, $L^{d,\infty}(\T^d)\subset L^s(\T^d)$ for all $s<d$. 

We shall use the endpoint pair $L^{d,\infty}(\T^d)$ and $L^{d/(d-1),1}(\T^d)$ in the following O'Neil--Young convolution inequality \cite{Blo72,ONe63}: There exists $C>0$ such that for all $f\in L^{d,\infty}(\T^d)$ and $g\in L^{d/(d-1),1}(\T^d)$,
\begin{align}\label{3.conv}
  \|f*g\|_{L^\infty(\T^d)} \le C\|f\|_{L^{d,\infty}(\T^d)}
  \|g\|_{L^{d/(d-1),1}(\T^d)}.
\end{align}
We also use the Lorentz--Sobolev inequality \cite{Alv77,SVS19}:
\begin{align}\label{3.lsi}
  \|f-\langle f\rangle\|_{L^{d/(d-1),1}(\T^d)}
  \le C\|\na f\|_{L^1(\T^d)}
\end{align}
for functions $f\in W^{1,1}(\T^d)$. The Zygmund space $L\log L(\T^d)$ consists of all measurable functions $f$ such that $\|f\|_{L\log L(\T^d)}<\infty$, where
\begin{align*}
  \|f\|_{L\log L(\T^d)} = \int_0^Mf^*(t)\log\frac{M}{t}\dd t,
  \quad M=|\T^d|.
\end{align*}
With this norm, $L\log L(\T^d)$ is a Banach space. For a nonnegative function $u$, it holds that $u\in L\log L(\T^d)$ if and only if $\int_{\T^d}u\log u\dx<\infty$. We denote by $\mbox{BMO}(\T^d)$ the space of functions of bounded mean oscillation, with its seminorm
\begin{align*}
  \|g\|_{\text{BMO}(\T^d)} = \sup_{Q\subset\T^d}\frac{1}{|Q|}
  \int_Q|g-\langle g\rangle_Q|\dx,
\end{align*}
where the supremum is taken over all cubes $Q\subset\T^d$ and 
$\langle g\rangle_Q$ is the mean value of $g$ over $Q$. We shall use the following $L\log L$--BMO pairing: There exists a constant $C>0$ such that for $f\in L\log L(\T^d)$ and $g\in{\rm BMO}(\T^d)$ with $\langle g\rangle=0$,
\begin{align}\label{3.pair}
  \bigg|\int_{\T^d}fg\dx\bigg| \le C\|f\|_{L\log L(\T^d)}
  \|g\|_{{\rm BMO}(\T^d)}.
\end{align}
Indeed, the dual space of $L\log L$ is $\exp L$ \cite[Chap.\ 4, Theorem 6.5]{BeRo88} and it holds that $\mbox{BMO}\hookrightarrow\exp L$ \cite{JoNi61}.

The classical Sobolev embedding $W^{1,p}(\T^d)\hookrightarrow L^\infty(\T^d)$ fails if $p=d$. Instead one may use the embedding into $\mbox{BMO}(\T^d)$, $\|f\|_{{\rm BMO}(\T^d)}\le C\|\na f\|_{L^d(\T^d)}$. This can be refined by the endpoint Lorentz space, as proven next.

\begin{lemma}[Endpoint Sobolev--BMO inequality]\label{lem.sbmo}
There exists a constant $C_d>0$ such that for all $g\in W^{1,1}_{\rm loc}(\T^d)$ satisfying $\na g\in L^{d,\infty}(\T^d;\R^d)$ and $\langle g\rangle=0$,
\begin{align*}
  \|g\|_{{\rm BMO}(\T^d)}\le C_d\|\na g\|_{L^{d,\infty}(\T^d)}.
\end{align*}
\end{lemma}

\begin{proof}
The result follows from the Poincar\'e inequality and the H\"older inequality for Lorentz spaces \cite[Chap.\ 4, Theorem 4.7]{BeRo88}. Indeed, for any cube $Q\subset\T^d$, we have
\begin{align*}
  \int_Q|g-\langle g\rangle_Q|\dx &\le C|Q|^{1/d}\int_Q|\na g|\dx
  \le C_d|Q|^{1/d}\|\na g\|_{L^{d,\infty}(\T^d)}
  \|{\rm 1}_Q\|_{L^{d/(d-1),1}(\T^d)} \\
  &\le C_d|Q|^{1/d}\|\na g\|_{L^{d,\infty}(\T^d)}|Q|^{(d-1)/d},
\end{align*}
implying that
\begin{align*}
  \frac{1}{|Q|}\int_Q|g-\langle g\rangle_Q|\dx
  \le C_d\|\na g\|_{L^{d,\infty}(\T^d)}.
\end{align*}
Taking the supremum over $Q$ yields the statement.
\end{proof}

%%%%%%%%%%%%%%%%%%%%

\subsection{Endpoint estimates}

We prove two endpoint estimates for the interaction energy and the drift field in terms of the $L^{d,\infty}(\T^d)$ norm of $\na K_{ij}$. 

\begin{lemma}[Endpoint control of the interaction energy]\label{lem.ec1}
Let Assumption \ref{ass.kern} hold. There exists a constant $C>0$ such that for all $u=(u_1,\ldots,u_n)$ with $u_i\in L\log L(\T^d)$ for $i=1,\ldots,n$, the interaction energy \eqref{eq:def-interaction-energy} is finite and satisfies
\begin{align*}
  |\mathcal P[u]|\le C\sum_{i,j=1}^n\|u_j\|_{L^1(\T^d)}
  \|\na K_{ij}\|_{L^{d,\infty}(\T^d)}\|u_i\|_{L\log L(\T^d)}.
\end{align*}
Moreover, if $u$ and $v$ are two vector-valued functions with finite entropy, then, for some constant $C>0$,
\begin{align*}
  |\mathcal P[u]-\mathcal P[v]| \le C\sum_{i,j=1}^n
  \|\na K_{ij}\|_{L^{d,\infty}(\T^d)}\|u_i+v_i\|_{L\log L(\T^d)}
  \|u_j-v_j\|_{L^1(\T^d)}. 
\end{align*}
\end{lemma}

\begin{proof}
Let $u_i\in L\log L(\T^d)$ for $i=1,\ldots,n$. By the $L\log L$--$\mbox{BMO}$ pairing \eqref{3.pair}, the BMO convolution inequality (which follows from Fubini's theorem), and Lemma \ref{lem.sbmo},
\begin{align*}
  \bigg|\int_{\T^d}u_i(K_{ij}*u_j)\dx\bigg|
  &\le C\|u_i\|_{L\log L(\T^d)}\|K_{ij}*u_j\|_{{\rm BMO}(\T^d)} \\
  &\le C\|u_i\|_{L\log L(\T^d)}\|K_{ij}\|_{{\rm BMO}(\T^d)}
  \|u_j\|_{L^1(\T^d)} \\
  &\le C\|u_i\|_{L\log L(\T^d)}\|\na K_{ij}\|_{L^{d,\infty}(\T^d)}
  \|u_j\|_{L^1(\T^d)}.
\end{align*}
The first statement results from summing these inequalities over $i,j=1,\ldots,n$.

We prove next the continuity estimate. It follows from the symmetry of $(K_{ij})$ as well as the $L\log L$--$\mbox{BMO}$ pairing \eqref{3.pair}, the BMO convolution inequality, and Lemma \ref{lem.sbmo} that
\begin{align*}
  |\mathcal P[u]-\mathcal P[v]|
  &= \bigg|\frac12\sum_{i,j=1}^n\int_{\T^d}(u_i+v_i)
  K_{ij}*(u_j-v_j)\dx\bigg| \\
  &\le C\sum_{i,j=1}^n\|u_i+v_i\|_{L\log L(\T^d)}
  \|K_{ij}*(u_j-v_j)\|_{{\rm BMO}(\T^d)} \\
  &\le C\sum_{i,j=1}^n\|u_i+v_i\|_{L\log L(\T^d)}
  \|\na K_{ij}\|_{L^{d,\infty}(\T^d)}\|u_j-v_j\|_{L^1(\T^d)}.
\end{align*}
This finishes the proof.
\end{proof}

\begin{lemma}[Endpoint control of the drift field]\label{lem.ec2}
Let Assumption \ref{ass.kern} hold. Then there exists a constant $C>0$ such that for all $u=(u_1,\ldots,u_n)\in W^{1,1}(\T^d;\R^n)$,
\begin{align*}
  \|\na p_i[u]\|_{L^\infty(\T^d)}
  \le C\sum_{j=1}^n\|\na K_{ij}\|_{L^{d,\infty}(\T^d)}
  \|\na u_j\|_{L^1(\T^d)}, \quad i=1,\ldots,n,
\end{align*}
recalling that $p_i[u] = \sum_{j=1}^n K_{ij}*u_j$.  
\end{lemma}

\begin{proof}
Since $\na K_{ij}$ is the weak gradient of a periodic function, it has zero mean, giving $\na K_{ij}*u_j = \na K_{ij}*(u_j-\langle u_j\rangle)$. Thus, using the O'Neil--Young convolution inequality \eqref{3.conv} and the Lorentz--Sobolev inequality \eqref{3.lsi},
\begin{align*}
  \|\na K_{ij}*u_j\|_{L^\infty(\T^d)}
  &\le C\|\na K_{ij}\|_{L^{d,\infty}(\T^d)}
  \|u_j-\langle u_j\rangle\|_{L^{d/(d-1),1}(\T^d)} \\
  &\le C\|\na K_{ij}\|_{L^{d,\infty}(\T^d)}
  \|\na u_j\|_{L^1(\T^d)},
\end{align*}
and a summation over $j=1,\ldots,n$ proves the lemma.
\end{proof}

Combining the previous lemma with the Fisher information bound $\|\na\sqrt{u_i}\|_{L^2(\T^d)}$ gives force and flux estimates.

\begin{corollary}[Force and flux bounds]\label{coro}
Let Assumption \ref{ass.kern} hold. There exists a constant $C>0$ such that for all $u=(u_1,\ldots,u_n)$ with finite mass $m_i=\int_{\T^d}u_i(t)\dx$ and $\na\sqrt{u_i}\in L^2((0,T)\times\T^d;\R^d)$, it holds that
\begin{align*}
  \|\na u_i\|_{L^2(0,T;L^1(\T^d))}
  &\le 2\sqrt{m_i}\|\na\sqrt{u_i}\|_{L^2((0,T)\times\T^d)}, \\
  \|\na p_i[u]\|_{L^2(0,T;L^\infty(\T^d))}
  &\le C\sum_{j=1}^n\|\na\sqrt{u_j}\|_{L^2((0,T)\times\T^d)}.
\end{align*}
In particular, the flux $F_i[u] = \na u_i+u_i\na p_i[u]$ satisfies
\begin{align*}
  \|F_i[u]\|_{L^2(0,T;L^1(\T^d))}
  \le C\sum_{j=1}^n\|\na\sqrt{u_j}\|_{L^2((0,T)\times\T^d)},  
\end{align*}
where the constant $C>0$ depends on the masses $m_j$ and $ \|\na K_{ij}\|_{L^{d,\infty}(\T^d)}$. 
\end{corollary}

\begin{proof}
The first estimate follows from $\na u_i=2\sqrt{u_i}\na\sqrt{u_i}$, while the second estimate is a consequence of the first one and Lemma \ref{lem.ec2}.
\end{proof}

%%%%%%%%%%%%%%%%%%%%%%%%%%%%%%%%%%%%%%%%%%%%%%%%%%%%%%%%%%%%%%

\section{Smooth approximation and existence of entropy--Fisher weak solutions}\label{sec.ex}

The existence of entropy--Fisher weak solutions was obtained in \cite{JLS26} using a finite volume approximation. We include here a proof based on regularized kernels, since this approximation is also used in the proof of the modulated stability estimate.

\subsection{Regularized approximation and entropy estimates}

We introduce a smooth approximation of the interaction kernel, which preserves the symmetry and positive semidefiniteness. Let $\rho$ be a standard even mollifier on $\R^d$ and let $\rho_\delta$ be the periodization of $\delta^{-d}\rho(\cdot/\delta)$ on $\T^d$ (not to be confused with the charge density $\rho[u]$), satisfying 
$\rho_\delta\in C^\infty(\T^d)$, $\rho_\delta\ge0$,
$\rho_\delta(x) = \rho_\delta(-x)$, and $\|\rho_\delta\|_{L^1(\T^d)}=1$.
Moreover, $\rho_\delta*f\to f$ in $L^p(\T^d)$ for $1\le p<\infty$, whenever $f\in L^p(\T^d)$. We define for $i,j=1,\ldots,n$:
\begin{align*}
  K_{ij}^\delta := \eta_\delta* K_{ij}, \quad
  \eta_\delta := \rho_\delta * \rho_\delta.
\end{align*}
Then $K_{ij}^\delta(x)=K_{ji}^\delta(-x)$, and the matrix kernel $K^\delta = (K_{ij}^\delta)_{i,j=1}^n$ remains positive semidefinite, since for every $f=(f_1,\dots,f_n)\in L^2(\T^d;\R^n)$,
\begin{align*}
  \sum_{i,j=1}^n\int_{\T^d}f_i(K_{ij}^\delta*f_j)\dx
  = \sum_{i,j=1}^n\int_{\T^d}(\rho_\delta*f_i)
  K_{ij}*(\rho_\delta*f_j)\dx \ge 0.
\end{align*}
We approximate the initial data by $u_{i,0}^\delta\in C^\infty(\T^d)$ such that $u_{i,0}^\delta$ is strictly positive, $\|u_{i,0}^\delta\|_{L^1(\T^d)} = \|u_{i,0}\|_{L^1(\T^d)} = m_i$, and
\begin{align}\label{4.init}
  \|u_{i,0}^\delta\|_{L\log L(\T^d)}\le C, \quad
  u_{i,0}^\delta\to u_{i,0}\quad\mbox{in }L^1(\T^d),
\end{align}
where the constant $C>0$ is independent of $\delta$.

For fixed $\delta>0$, we consider the regularized system
\begin{align}\label{4.approx}
  \pa_t u_i^\delta = \Delta u_i^\delta 
  + \diver(u_i^\delta\na p_i^\delta[u^\delta]), \quad
  p_i^\delta[u^\delta] = \sum_{j=1}^n K_{ij}^\delta * u_j^\delta,
\end{align}
with initial condition $u_i^\delta(0) = u_{i,0}^\delta$ in $\T^d$ for $i=1,\ldots,n$. Let $m=(m_1,\ldots,m_n)$ and introduce the relative entropy
\begin{align}\label{4.E}
  \mathcal E[u^\delta|m] = \sum_{i=1}^n\int_{\T^d}u_i^\delta
  \log\frac{u_i^\delta}{m_i}\dx
\end{align}
as well as the regularized interaction energy and dissipation, respectively,
\begin{align*}
  \mathcal P^\delta[u] = \frac12\sum_{i=1}^n\int_{\T^d}
  u_ip_i^\delta[u]\dx, \quad
  \mathcal D^\delta[u] = \sum_{i=1}^n\int_{\T^d}\na p_i^\delta[u]
  \cdot\na u_i \dx.
\end{align*}
Observe that the positive semidefiniteness of $K^\delta$ implies that $\mathcal P^\delta[u^\delta]\ge 0$ and $\mathcal D^\delta[u^\delta]\ge 0$. We prove now the existence of a solution to the regularized problem. 

\begin{proposition}[Regularized solutions and uniform estimates]
\label{prop.delta}
System \eqref{4.approx} admits a global positive strong solution $u^\delta= (u_1^\delta,\ldots,u_n^\delta)\in W^{1,\infty}((0,T)\times\T^d;\R^n)$. For every $t\ge 0$, the masses are conserved, $\|u_i^\delta(t)\|_{L^1(\T^d)}=m_i$ for $i=1,\ldots,n$, and the relative entropy \eqref{4.E} satisfies
\begin{align}\label{4.delta}
  \mathcal E[u^\delta(t)|m] + 4\sum_{i=1}^n\int_0^t\int_{\T^d}
  |\na(u_i^\delta)^{1/2}|^2\dx \dd s
  + \int_0^t\mathcal D^\delta[u^\delta(s)]\dd s
  = \mathcal E[u_0^\delta|m].
\end{align}
\end{proposition}

We deduce from the relative entropy identity and Corollary \ref{coro} that there exists a constant $C>0$ such that the following uniform estimates hold:
\begin{align}\label{4.est}
  \sup_{0<t<T}&\sum_{i=1}^n\|u_i^\delta(t)\|_{L\log L(\T^d)}
  + \sum_{i=1}^n\|\na u_i^\delta\|_{L^2(0,T;L^1(\T^d))} \\
  &+ \sum_{i=1}^n\|\na p_i^\delta[u^\delta]\|_{L^2(0,T;L^\infty(\T^d))}
  + \sum_{i=1}^n\|F_i^\delta[u^\delta]\|_{L^2(0,T;L^1(\T^d))} \le C,
  \nonumber 
\end{align}
where $F_i^\delta[u^\delta] = \na u_i^\delta + u_i^\delta\na p_i^\delta[u^\delta]$ is the flux of the $i$th species. 

\begin{proof}[Proof of Proposition \ref{prop.delta}]
Since the kernels $K_{ij}^\delta$ and the initial data $u_{i,0}^\delta$ are smooth, standard semilinear parabolic theory gives a unique local smooth solution. The divergence form implies that the masses are conserved. Setting $B_i^\delta = \na p_i^\delta[u^\delta]$ for the given solution $u^\delta$, it solves the scalar linear equation
\begin{align*}
  \pa_t u_i^\delta = \Delta u_i^\delta + B_i^\delta\cdot\na u_i^\delta
  + c_i^\delta u_i^\delta, \quad\mbox{where } 
  c_i^\delta = \diver B_i^\delta.
\end{align*}
On every finite time interval, the coefficients are smooth and $c_i^\delta\le M$ is bounded. After the change of unknowns $w_i^\delta = e^{-Mt}u_i^\delta$, since $u_{i,0}^\delta>0$ by construction, we can apply the strong maximum principle to $w_i^\delta$ and obtain $u_i^\delta(t,x)>0$ for all $t>0$, $x\in\T^d$. 

The local solution can be extended to a global one. Indeed, the weak maximum-principle argument gives
\begin{align*}
  \|u_i^\delta(t)\|_{L^\infty(\T^d)} \le e^{C(\delta)t}
  \|u_{i,0}^\delta\|_{L^\infty(\T^d)}.
\end{align*}
Noting that the regularized kernels are smooth, we can also derive bounds for $\|u_i^\delta(t)\|_{W^{1,\infty}(\T^d)}$. These bounds prevent the finite-time blow-up by the parabolic continuation criterion. Hence, the smooth solution exists globally for $t>0$. Finally, the relative entropy identity is obtained after multiplying \eqref{4.approx} by $\log(u_i^\delta/m_i)$, integrating over $\T^d$, and summing over $i=1,\ldots,n$:
\begin{align}\label{4.dEdt}
  \frac{\dd}{\dd t}\mathcal E[u^\delta|m]
  + \sum_{i=1}^n\int_{\T^d}\frac{|\na u_i^\delta|^2}{u_i^\delta}\dx
  + \sum_{i=1}^n\int_{\T^d}\na u_i^\delta\cdot\na p_i^\delta[u^\delta]
  \dx = 0,
\end{align}
which gives \eqref{4.delta} after integration over $(0,t)$. 
\end{proof}

%%%%%%%%%%%%%%%%%

\subsection{Compactness and existence of entropy--Fisher weak solutions}

We pass to the limit $\delta\to 0$ in the regularized system \eqref{4.approx}. Proposition \ref{prop.ex} is a direct consequence of the following result.

\begin{theorem}[Existence of entropy--Fisher weak solutions]
\label{thm.ex} 
Let Assumptions \ref{ass.kern} and \ref{ass.init} hold and let $u^\delta$ be the regularized strong solution to \eqref{4.approx} constructed in Proposition \ref{prop.delta}. Then, for any $T>0$, there exists a subsequence of $u^\delta$ (not relabeled) such that, as $\delta\to 0$,
\begin{align*}
  u_i^\delta\to u_i\quad\mbox{strongly in }L^2(0,T;L^1(\T^d)),
\end{align*}
and $u=(u_1,\ldots,u_n)$ is an entropy--Fisher weak solution to \eqref{1.eq}--\eqref{1.ic}.
\end{theorem}

\begin{proof}
The uniform estimates \eqref{4.est} show that $(u_i^\delta)$ is bounded in $L^2(0,T;W^{1,1}(\T^d))$ and $(F_i^\delta[u^\delta])$ is bounded in $L^2(0,T;L^1(\T^d))$. Consequently, $\pa_t u_i^\delta = \diver F_i^\delta[u^\delta]$ is uniformly bounded in $L^2(0,T;W^{-1,1}(\T^d))$. Since the embedding $W^{1,1}(\T^d)\hookrightarrow L^1(\T^d)$ is compact, the Aubin--Lions lemma implies the existence of a subsequence, which is not relabeled, such that, as $\delta\to 0$,
\begin{align*}
  u_i^\delta \to u_i\quad\mbox{strongly in }L^2(0,T;L^1(\T^d)). 
\end{align*}
We deduce from $u_i^\delta>0$ that $u_i\ge 0$ a.e.\ in $(0,T)\times\T^d$. 

To identify the limiting diffusive flux, we use $(\sqrt{a}-\sqrt{b})^2\le |a-b|$ to infer that $(u_i^\delta)^{1/2}\to u_i^{1/2}$ strongly in $L^2((0,T)\times\T^d)$. The uniform Fisher-information bound implies that, up to a subsequence,
\begin{align}\label{4.H1}
  (u_i^\delta)^{1/2}\rightharpoonup u_i^{1/2}
  \quad\mbox{weakly in }L^2(0,T;H^1(\T^d)). 
\end{align}
We conclude that
\begin{align}\label{4.nau}
  \na u_i^\delta = 2(u_i^\delta)^{1/2}\na(u_i^\delta)^{1/2}
  \rightharpoonup 2u_i^{1/2}\na u_i^{1/2} = \na u_i
  \quad\mbox{weakly in }L^1((0,T)\times\T^d).
\end{align}

We next identify the drift fields. We know from \eqref{4.est} that $(\na p_i^\delta[u^\delta])$ is bounded in $L^2(0,T;L^\infty(\T^d;\R^d))$. Thus, up to a subsequence,
\begin{align*}
  \na p_i^\delta[u^\delta]\rightharpoonup \zeta_i\quad
  \mbox{weakly* in }L^2(0,T;L^\infty(\T^d;\R^d)).
\end{align*}
We claim that $\zeta_i=\na p_i[u]$. Indeed, we write
\begin{align}\label{4.aux}
  \na p_i^\delta[u^\delta] - \na p_i[u]
  = (\na p_i^\delta[u^\delta] - \na p_i^\delta[u])
  + (\na p_i^\delta[u] - \na p_i[u]).
\end{align}
Consider the second term. Since $u_i\in L^2(0,T;W^{1,1}(\T^d))$, the family $(\eta_\delta*u_i)$ converges to $u_i$ strongly in $L^2(0,T;W^{1,1}(\T^d))$. Hence, by Lemma \ref{lem.ec2},
\begin{align*}
  \na p_i^\delta[u] = \na p_i[\eta_\delta * u]\to \na p_i[u]
  \quad\mbox{strongly in }L^2(0,T;L^\infty(\T^d)). 
\end{align*}
For the first term in \eqref{4.aux}, we use the evenness of the mollifier $\rho_\delta$ and the symmetry of $K_{ij}^\delta$, and move the convolution to the test function $\Phi\in C^\infty((0,T)\times\T^d;\R^d)$:
\begin{align}\label{4.aux2}
  \int_0^T\int_{\T^d}\big(\na p_i^\delta[u^\delta] 
  - \na p_i^\delta[u]\big)\cdot\Phi\dx\dd t
  = -\sum_{j=1}^n\int_0^T\int_{\T^d}(u_j^\delta-u_j)
  (\na K_{ji}^\delta*\Phi)\dx\dd t.
\end{align}
The term $\na K_{ji}^\delta*\Phi = \na K_{ji}*(\eta_\delta*\Phi)$ is uniformly bounded in $L^\infty((0,T)\times\T^d)$, by Lemma \ref{lem.ec2}. Therefore, in view of $u_j^\delta\to u_j$ strongly in $L^1((0,T)\times\T^d)$, the right-hand side of \eqref{4.aux2} vanishes in the limit $\delta\to 0$. This proves that 
\begin{align*}
  \na p_i^\delta[u^\delta]\to \na p_i[u] \quad\mbox{in the sense
  of distributions},
\end{align*}
and consequently $\zeta_i=\na p_i[u]$. The strong convergence of $u_i^\delta$ and the weak convergence of $\na p_i^\delta[u^\delta]$ imply that
\begin{align*}
  u_i^\delta\na p_i^\delta[u^\delta]\rightharpoonup u_i\na p_i[u]
  \quad\mbox{weakly in }L^1((0,T)\times\T^d). 
\end{align*}
Taking into account convergence \eqref{4.nau}, we have $F_i^\delta[u^\delta]\rightharpoonup F_i[u]$ weakly in $L^1((0,T)\times\T^d)$. Therefore, passing to the limit in the weak formulation of \eqref{4.approx}, for every $\psi\in C^\infty([0,T]\times\T^d)$ with $\psi(T)=0$,
\begin{align*}
  -\int_0^T\int_{\T^d}u_i\pa_t\psi\dx\dd t 
  - \int_{\T^d}u_{i,0}\psi(0)\dx 
  = -\int_0^T\int_{\T^d}F_i[u]\cdot\na\psi\dx\dd t.
\end{align*}
Hence, the limit $u$ is a distributional solution.

It remains to verify that the limit is an entropy--Fisher weak solution. Since $u_i\ge 0$ and $u_i^\delta(t)\to u_i(t)$ strongly in $L^1(\T^d)$ for a.e.\ $t>0$, the lower semicontinuity of the relative entropy gives
\begin{align*}
  \mathcal E[u(t)|m] \le \liminf_{\delta\to 0}\mathcal E
  [u^\delta(t)|m] \quad\mbox{for a.e. }t\in(0,T).
\end{align*}
Moreover, using convergence \eqref{4.H1} and the weak lower semicontinuity of the Fisher information,
\begin{align*}
  \sum_{i=1}^n\int_0^t\int_{\T^d}|\na u_i^{1/2}|^2\dx\dd s
  \le \liminf_{\delta\to 0}\sum_{i=1}^n\int_0^t\int_{\T^d}
  |\na(u_i^\delta)^{1/2}|^2\dx\dd s.
\end{align*}
Then the entropy--Fisher bound follows from the entropy identity in Proposition \ref{prop.delta}. By Corollary \ref{coro}, we have the regularity $\na p_i[u]\in L^2(0,T;L^\infty(\T^d;\R^d))$ and $F_i[u]\in L^2(0,T;$ $L^1(\T^d;\R^d))$. Thus, the weak formulation extends to all test functions $\psi\in W^{1,\infty}((0,T)\times\T^d)$ with $\psi(T)=0$. This finishes the proof.
\end{proof}

%%%%%%%%%%%%%%%%%%%%%%%%%%%%%%%%%%%%%%%%%%%%%%%%%%%%%%%%%%%%%%

\section{Renormalized Hellinger distance}\label{sec.hell}

The aim of this section is the estimation of the Hellinger distance \eqref{1.HP}. First, we derive a renormalized formulation of \eqref{1.eq}.

\subsection{Renormalized formulation}

To derive a renormalized formulation for \eqref{1.eq}, we apply a similar strategy as in \cite{Fis17}. This formulation is used later to justify the nonlinear entropy identities involving a regularized comparison solution. 

\begin{lemma}[Renormalized formulation of entropy--Fisher weak solutions]
\label{lem.renorm}
Let $u=(u_1,\ldots,$ $u_n)$ be an entropy--Fisher weak solution to \eqref{1.eq} on $(0,T)\times\T^d$ and let $\xi\in C^2([0,\infty))$ satisfy $\xi'\in C^1_c([0,\infty))$. Then, for all $i=1,\ldots,n$, $\varphi\in W^{1,\infty}([0,T]\times\T^d)$, and for a.e.\ $\tau\in(0,T]$, $u$ satisfies
\begin{align}\label{5.xi}
  \int_{\T^d}&\xi(u_i(\tau))\varphi(\tau)\dx
  - \int_{\T^d}\xi(u_{i,0})\varphi(0)\dx
  - \int_0^\tau\int_{\T^d}\xi(u_i)\pa_t\varphi\dx\dd t \\
  &= -\int_0^\tau\int_{\T^d}\xi'(u_i)F_i[u]\cdot\na\varphi\dx\dd t
  - \int_0^\tau\int_{\T^d}\xi''(u_i)F_i[u]\cdot\na u_i\varphi\dx\dd t,
  \nonumber 
\end{align}
recalling that $F_i[u]=\na u_i+u_i\na p_i[u]$. 
\end{lemma}

\begin{proof}
Let $\rho_\eta$ be a standard nonnegative spatial mollifier on $\T^d$ and set $u_i^\eta:=\rho_\eta*u_i$ and $F_i^\eta:=\rho_\eta*F_i[u]$. The weak formulation of \eqref{1.eq} for $u_i^\eta$ reads as $\pa_t u_i^\eta = \diver F_i^\eta$. Since $u_i^\eta$ and $F_i^\eta$ are smooth in space, we find for any $\xi$ that
\begin{align*}
  \pa_t\xi(u_i^\eta) = \diver\big(\xi'(u_i^\eta)F_i^\eta\big)
  - \xi''(u_i^\eta)F_i^\eta\cdot\na u_i^\eta. 
\end{align*}
We multiply this identity by $\varphi\in W^{1,\infty}([0,T]\times\T^d)$ and integrate over $(0,\tau)\times\T^d$ for some $\tau>0$:
\begin{align}\label{5.aux}
  \int_{\T^d}&\xi(u_i^\eta(\tau))\varphi(\tau)\dx
  - \int_{\T^d}\xi(u_i^\eta(0))\varphi(0)\dx
  - \int_0^\tau\int_{\T^d}\xi(u_i^\eta)\pa_t\varphi\dx\dd t \\
  &= -\int_0^\tau\int_{\T^d}\xi'(u_i^\eta)F_i^\eta\cdot\na\varphi\dx\dd t
  - \int_0^\tau\int_{\T^d}\xi''(u_i^\eta)F_i^\eta\cdot\na u_i^\eta
  \varphi\dx\dd t. \nonumber 
\end{align}
We pass to the limit $\eta\to 0$. Since $u_i^\eta(\tau)\to u_i(\tau)$ in $L^1(\T^d)$ for a.e.\ $\tau\in(0,T)$, $\xi$ is Lipschitz continuous, and $\varphi$ is bounded, the left-hand side of \eqref{5.aux} converges to the left-hand side of \eqref{5.xi}. For the first term on the right-hand side, we have
\begin{align*}
  \int_0^T\int_{\T^d}|\xi'(u_i^\eta)F_i^\eta - \xi'(u_i)F_i[u]|\dx\dd t
  &\le \int_0^T\int_{\T^d}|\xi'(u_i^\eta)||F_i^\eta-F_i|\dx\dd t \\
  &\phantom{xx}
  + \int_0^T\int_{\T^d}|\xi'(u_i^\eta)-\xi'(u_i)||F_i|\dx\dd t.
\end{align*}
As $\xi'(u_i^\eta)$ is uniformly bounded and $F_i^\eta$ converges to $F_i$ strongly in $L^1((0,T)\times\T^d)$, the first term on the right-hand side converges to zero as $\eta\to 0$. 
Moreover, the strong convergence of $u_i^\eta$ in $L^1((0,T)\times\T^d)$ implies that $u_i^\eta\to u_i$ in measure. Hence $\xi'(u_i^\eta)$ converges to $\xi'(u_i)$ in measure by the Lipschitz continuity of $\xi'$. Note that $|\xi'(u_i^\eta)F_i[u]|\le C|F_i[u]|$, i.e., it is bounded by an $L^1((0,T)\times\T^d)$ function. Vitali's version of dominated convergence shows that the second term on the right-hand side converges to zero. We have proved that $\xi'(u_i^\eta)F_i^\eta\to\xi'(u_i)F_i[u]$ strongly in $L^1((0,T)\times\T^d)$. Hence, the first term on the right-hand side of \eqref{5.aux} converges. 

It remains to pass to the limit in the second term of \eqref{5.aux}. Note that $|F_i^\eta|\le |\na u_i^\eta|+ |\rho_\eta*(u_i\na p_i[u])|$. We estimate, using H\"older's inequality,
\begin{align*}
  |\na u_i^\eta| = 2|\rho_\eta*(\sqrt{u_i}\na\sqrt{u_i})|
  &\le 2(u_i^\eta)^{1/2}\big(\rho_\eta*|\na\sqrt{u_i}|^2\big)^{1/2}, \\
  |\rho_\eta*(u_i\na p_i[u])|
  &\le (u_i^\eta)^{1/2}\big(\rho_\eta*(u_i|\na p_i[u]|^2)\big)^{1/2}.
\end{align*}
Let $R>0$ be such that $\xi'\in C_c^1([0,R))$. Since $\xi''(u_i^\eta)=0$ outside $\{u_i^\eta\le R\}$, we infer that
\begin{align*}
  |\xi''(u_i^\eta)F_i^\eta\cdot\na u_i^\eta|
  \le C_R\big(\rho_\eta*|\na\sqrt{u_i}|^2 
  + \rho_\eta*(u_i|\na p_i[u]|^2)\big).
\end{align*}
The right-hand side is uniformly integrable on $(0,T)\times\T^d$ and $u_i^\eta\to u_i$, $F_i^\eta\to F_i[u]$, $\na u_i^\eta\to \na u_i$ in measure on $(0,T)\times\T^d$. So, we can argue as above to conclude that
\begin{align*}
  \xi''(u_i^\eta)F_i^\eta\cdot\na u_i^\eta
  \to \xi''(u_i)F_i[u]\cdot\na u_i
  \quad\mbox{strongly in }L^1((0,T)\times\T^d).
\end{align*}
Hence, the last term in \eqref{5.aux} also converges. This finishes the proof.
\end{proof}

%%%%%%%%%%%%%%%%%%

\subsection{Square-root renormalization}

The formulation in Lemma \ref{lem.renorm} is stated for renormalized functions whose derivatives are compactly supported. We extend it to the square-root renormalization. 

\begin{lemma}[Square-root renormalization]\label{lem.renorm2}
Let $u=(u_1,\ldots,u_n)$ be an entropy--Fisher weak solution to \eqref{1.eq}. Set $\eta(s)=\sqrt{s+\eps}$ for $s\ge 0$, $\eps>0$. Then, for all $i=1,\ldots,n$, $\varphi\in W^{1,\infty}((0,T)\times\T^d)$ and for a.e.\ $\tau\in(0,T]$, $u$ satisfies
\begin{align*}
  \int_{\T^d}&\eta(u_i(\tau))\varphi(\tau)\dx
  - \int_{\T^d}\eta(u_{i,0})\varphi(0)\dx
  - \int_0^\tau\int_{\T^d}\eta(u_i)\pa_t\varphi\dx\dd t \\
  &= -\int_0^\tau\int_{\T^d}\eta'(u_i)F_i[u]\cdot\na\varphi\dx\dd t
  - \int_0^\tau\int_{\T^d}\eta''(u_i)F_i[u]\cdot\na u_i\varphi\dx\dd t.
  \nonumber 
\end{align*}
\end{lemma}

\begin{proof}
Let $\chi_M\in C_c^\infty([0,\infty);[0,1])$ for $M>0$ be such that $0\le\chi_M\le 1$ and $\chi_M=1$ on $[0,M]$, $\chi_M=0$ on $[2M,\infty)$, and $|\chi'_M|\le C/M$. We approximate $\eta$ by setting
\begin{align*}
  \eta_M(s) = \eta(0) + \int_0^s\chi_M(r)\eta'(r)\dd r
  \quad\mbox{for }s\ge 0.
\end{align*}
Then $\eta'_M\in C_c^1([0,\infty))$, and the renormalized identity \eqref{5.xi} holds for $\xi=\eta_M$. We wish to pass to the limit $M\to\infty$. 

It holds that $\eta_M(s)\to\eta(s)$ and $\eta'_M(s)=\chi_M(s)\eta'(s)\to \eta'(s)$ pointwise for every $s\ge 0$. Since $F_i[u]\in L^1((0,T)\times\T^d;\R^d)$ and $\eta(u_i)\in L^1((0,T)\times\T^d)$, together with the bounds
\begin{align*}
  0\le\eta_M(s)\le\eta(s), \quad |\eta_M'(s)|\le C\eps^{-1/2},
\end{align*}
the terms involving $\eta_M(u_i)$ and $\eta'_M(u_i)$ in \eqref{5.xi} converge by dominated convergence. It remains to treat the term involving $\xi''=\eta_M''$. We have $\eta''_M=\chi_M\eta''+\chi_M'\eta'$ and 
\begin{align*}
  |\eta''(u_i)F_i[u]\cdot\na u_i|
  &\le |\eta''(u_i)|\big(|\na u_i|^2 + |u_i\na p_i[u]||\na u_i|\big) \\
  &\le C_\eps|\na\sqrt{u_i}|^2 
  + C_\eps|\sqrt{u_i}\na p_i[u]||\na\sqrt{u_i}|.
\end{align*}
The right-hand side is integrable on $(0,T)\times\T^d$ by the definition of the entropy--Fisher weak solution and Corollary \ref{coro}. Hence, 
\begin{align*}
  \chi_M(u_i)\eta''(u_i)F_i[u]\cdot\na u_i
  \to \eta''(u_i)F_i[u]\cdot\na u_i
\end{align*}
strongly in $L^1((0,T)\times\T^d)$ as $M\to\infty$. The other cutoff contribution $\chi_M'\eta'$ is supported on $\{M\le u_i\le 2M\}$. On this set, we have $|\chi'_M(u_i)\eta'(u_i)|\le C/(M\sqrt{u_i+\eps})$. We obtain
\begin{align*}
  |\chi'_M(u_i)\eta'(u_i)||F_i[u]||\na u_i|
  \le \frac{C}{\sqrt{M}}\big(|\na\sqrt{u_i}|^2
  + |\sqrt{u_i}\na p_i[u]||\na\sqrt{u_i}|\big).
\end{align*}
The right-hand side belongs to $L^1((0,T)\times\T^d)$, by Corollary \ref{coro} again. Hence, this contribution converges to zero strongly in $L^1((0,T)\times\T^d)$. Thus, the limit $M\to\infty$ in the renormalized formulation for $\eta_M$ completes the proof.
\end{proof}

%%%%%%%%%%%%%%%%%%%

\subsection{Regularized Hellinger identity}

We derive the Hellinger identity. To this end, we compare an entropy--Fisher weak solution $u=(u_1,\ldots,u_n)$ to \eqref{1.eq} with a regularized solution $v^\delta=(v_1^\delta,\ldots,v_n^\delta)$ to \eqref{4.approx}. We introduce the functions
\begin{align*}
  U_i^\eps = u_i+\eps, \quad V_i^\eps = v_i^\delta+\eps, \quad
  L_i^\eps = \log\frac{U_i^\eps}{V_i^\eps}\quad\mbox{for }0<\eps<1.
\end{align*}
Since entropy--Fisher weak solutions may vanish, we define the regularized Hellinger distance
\begin{align*}
  \mathcal H_\eps[u|v^\delta] = \mathcal H[u+\eps|v^\delta+\eps]
  = \sum_{i=1}^n\int_{\T^d}\big(\sqrt{U_i^\eps}
  - \sqrt{V_i^\eps}\big)^2\dx. 
\end{align*}

\begin{proposition}[Regularized Hellinger identity]\label{prop.rhi}
Let $u$ be an entropy--Fisher weak solution to \eqref{1.eq} and $v^\delta$ be a regularized solution in the sense of Proposition \ref{prop.delta}. Then, for a.e.\ $\tau\in (0,T]$, 
\begin{align}\label{5.rhi}
  \mathcal H_\eps&[u|v^\delta](\tau)
  - \mathcal H_\eps[u_0|v_0^\delta]
  + \frac12\sum_{i=1}^n\int_0^\tau\int_{\T^d}
  \sqrt{U_i^\eps V_i^\eps}|\na L_i^\eps|^2\dx\dd t \\
  &= -\frac12\sum_{i=1}^n\int_0^\tau\int_{\T^d}
  \sqrt{U_i^\eps V_i^\eps}\na(p_i[u]-p_i^\delta[v^\delta])
  \cdot\na L_i^\eps\dx\dd t \nonumber \\
  &\phantom{xx}+ \frac{\eps}{2}\sum_{i=1}^n\int_0^\tau\int_{\T^d}
  \sqrt{U_i^\eps V_i^\eps}\bigg(\frac{\na p_i[u]}{U_i^\eps}
  -\frac{\na p_i^\delta[v^\delta]}{V_i^\eps}\bigg)\cdot\na L_i^\eps 
  \dx\dd t. \nonumber 
\end{align}
\end{proposition}

\begin{proof}
We apply Lemma \ref{lem.renorm2} by choosing the test function $\varphi=\sqrt{V_i^\eps}$. This test function is admissible since $\sqrt{V_i^\eps}\in W^{1,\infty}((0,T)\times\T^d)$. It follows from the identity $\widetilde\eta'(u_i)\na\varphi + \widetilde\eta''(u_i)\na u_i\varphi = \na(\widetilde\eta'(u_i)\varphi)$, where $\widetilde\eta(s)=\sqrt{s+\eps}$, that
\begin{align*}%\label{5.aux2}
  \int_{\T^d}&\sqrt{U_i^\eps V_i^\eps}(\tau)\dx
  - \int_{\T^d}\sqrt{U_i^\eps V_i^\eps}(0)\dx \\
  &= -\frac12\int_0^\tau\int_{\T^d}F_i[u]\cdot\na
  \bigg(\frac{V_i^\eps}{U_i^\eps}\bigg)^{1/2}\dx\dd t
  + \int_0^\tau\int_{\T^d}\sqrt{U_i^\eps}\pa_t\sqrt{V_i^\eps}\dx\dd t.
  \nonumber 
\end{align*}
Replacing $\pa_t\sqrt{V_i^\eps}$ in the last term by
\begin{align*}
  \pa_t\sqrt{V_i^\eps} 
  = \frac{\diver F_i^\delta[v^\delta]}{2\sqrt{V_i^\eps}}, \quad
  \mbox{where }F_i^\delta[v^\delta] = \na v_i^\delta + v_i^\delta
  \na p_i^\delta[v^\delta],
\end{align*}
summing over $i=1,\ldots,n$, using
\begin{align*}
  \mathcal H_\eps[u|v^\delta] = \sum_{i=1}^n\int_{\T^d}
  (u_i+v_i^\delta+2\eps)\dx - 2\sum_{i=1}^n\int_{\T^d}
  \sqrt{U_i^\eps V_i^\eps}\dx,
\end{align*}
and observing that both $u_i$ and $v_i^\delta$ conserve the mass, we find that
\begin{align}\label{5.aux2}
  \mathcal H_\eps[u|v^\delta](\tau) - \mathcal H_\eps[u_0|v_0^\delta]
  &= \sum_{i=1}^n\int_0^\tau\int_{\T^d}F_i[u]\cdot\na
  \bigg(\frac{V_i^\eps}{U_i^\eps}\bigg)^{1/2}\dx\dd t \\
  &\phantom{xx}+ \sum_{i=1}^n\int_0^\tau\int_{\T^d}F_i^\delta[v^\delta]\cdot\na
  \bigg(\frac{U_i^\eps}{V_i^\eps}\bigg)^{1/2}\dx\dd t. \nonumber 
\end{align}
We compute the two terms on the right-hand side. The diffusion terms $\na u_i=\na U_i^\eps$ and $\na v_i^\delta=\na V_i^\eps$ combine into
\begin{align*}
  \na U_i^\eps\cdot\na\bigg(\frac{V_i^\eps}{U_i^\eps}\bigg)^{1/2}
  + \na V_i^\eps\cdot\na\bigg(\frac{U_i^\eps}{V_i^\eps}\bigg)^{1/2}
  = -\frac12\sqrt{U_i^\eps V_i^\eps}|\na L_i^\eps|^2.
\end{align*}
Furthermore, the drift terms become
\begin{align*}
  u_i&\na p_i[u]\cdot\na\bigg(\frac{V_i^\eps}{U_i^\eps}\bigg)^{1/2}
  + v_i^\delta\na p_i^\delta[v^\delta]\cdot  \na\bigg(\frac{U_i^\eps}{V_i^\eps}\bigg)^{1/2} \\
  &= -\frac12\bigg\{\bigg(\frac{V_i^\eps}{U_i^\eps}\bigg)^{1/2}
  u_i\na p_i[u] - \bigg(\frac{U_i^\eps}{V_i^\eps}\bigg)^{1/2}
  v_i^\delta\na p_i^\delta[v^\delta]\bigg\}\cdot\na L_i^\eps \\
  &= -\frac12\sqrt{U_i^\eps V_i^\eps}\na(p_i[u]-p_i^\delta[v^\delta])
  \cdot\na L_i^\eps + \frac12\eps\sqrt{U_i^\eps V_i^\eps}
  \bigg(\frac{\na p_i[u]}{U_i^\eps}
  - \frac{\na p_i^\delta[v^\delta]}{V_i^\eps}\bigg)
  \cdot\nabla L_i^\varepsilon.
\end{align*}
Substituting these identities into \eqref{5.aux2} yields \eqref{5.rhi} and concludes the proof.
\end{proof}

%%%%%%%%%%%%%%%%%%%%%

\subsection{Estimated Hellinger identity}

We estimate the regularized Hellinger identity derived in Proposition \ref{prop.rhi}. 

\begin{lemma}[Estimated Hellinger identity]\label{lem.ehi}
Let $u$ be an entropy--Fisher weak solution to \eqref{1.eq} and $v^\delta$ be a regularized solution in the sense of Proposition \ref{prop.delta}. Then there exists a constant $C>0$, independent of $\eps$, such that for a.e.\ $\tau\in(0,T]$, 
\begin{align}\label{5.ehi}
  \mathcal H_\eps&[u|v^\delta](\tau)
  - \mathcal H_\eps[u_0|v_0^\delta]
  + \frac14\sum_{i=1}^n\int_0^\tau\int_{\T^d}
  \sqrt{U_i^\eps V_i^\eps}|\na L_i^\eps|^2\dx\dd t \\
  &\le \frac12\sum_{i=1}^n\int_0^\tau\int_{\T^d}
  \sqrt{U_i^\eps V_i^\eps}|\na(p_i[u]-p_i^\delta[v^\delta])|^2
  \dx\dd t \nonumber \\
  &\phantom{xx} + C\sqrt{\eps}\sum_{i=1}^n\int_0^\tau
  \big(\|\na p_i[u]\|_{L^\infty(\T^d)}^2
  + \|\na p_i^\delta[v^\delta]\|_{L^\infty(\T^d)}^2\big)\dd t.
  \nonumber 
\end{align}
\end{lemma}

\begin{proof}
We estimate the two terms on the right-hand side of \eqref{5.rhi}. By Young's inequality, the first one becomes
\begin{align*}
  \bigg|\frac12&\sum_{i=1}^n\int_0^\tau\int_{\T^d}
  \sqrt{U_i^\eps V_i^\eps}\na(p_i[u]-p_i^\delta[v^\delta])
  \cdot\na L_i^\eps\dx\dd t\bigg| \\
  &\le \frac18\sum_{i=1}^n\int_0^\tau\int_{\T^d}\sqrt{U_i^\eps V_i^\eps}
  |\na L_i^\eps|^2\dx\dd t
  + \frac12 \sum_{i=1}^n\int_0^\tau\int_{\T^d}\sqrt{U_i^\eps V_i^\eps}
  |\na(p_i[u]-p_i^\delta[v^\delta])|^2\dx\dd t,
\end{align*}
and the first term on the right-hand side can be absorbed by the left-hand side of \eqref{5.rhi}. Similarly, the other term is estimated as
\begin{align*}
  \bigg|\frac{\eps}{2}&\sum_{i=1}^n\int_0^\tau\int_{\T^d}
  \sqrt{U_i^\eps V_i^\eps}\bigg(\frac{\na p_i[u]}{U_i^\eps}
  -\frac{\na p_i^\delta[v^\delta]}{V_i^\eps}\bigg)\cdot\na L_i^\eps 
  \dx\dd t\bigg| \\
  &\le \frac18\sum_{i=1}^n\int_0^\tau\int_{\T^d}\sqrt{U_i^\eps V_i^\eps}
  |\na L_i^\eps|^2\dx\dd t + R_\eps(\tau),
\end{align*}
where
\begin{align*}
  R_\eps(\tau) := \frac{\eps^2}{2}\sum_{i=1}^n\int_0^\tau\int_{\T^d}
  \sqrt{U_i^\eps V_i^\eps}\bigg|\frac{\na p_i[u]}{U_i^\eps}
  -\frac{\na p_i^\delta[v^\delta]}{V_i^\eps}\bigg|^2\dx\dd t.
\end{align*}
It remains to estimate $R_\eps(\tau)$. Because of $U_i^\eps$, $V_i^\eps\ge \eps$, we have
\begin{align*}
  \eps^2\frac{\sqrt{U_i^\eps V_i^\eps}}{(U_i^\eps)^2}
  \le \sqrt{\eps}\sqrt{V_i^\eps} \le C\sqrt{\eps}(1+V_i^\eps), \quad
  \eps^2\frac{\sqrt{U_i^\eps V_i^\eps}}{(V_i^\eps)^2}
  \le C\sqrt{\eps}(1+U_i^\eps).
\end{align*}
Therefore, since the $L^1(\T^d)$ norms of $U_i^\eps$ and $V_i^\eps$ are uniformly bounded,
\begin{align*}
  R_\eps(\tau) \le C\sqrt{\eps}\sum_{i=1}^n\int_0^\tau
  \big(\|\na p_i[u]\|_{L^\infty(\T^d)}^2
  + \|\na p_i^\delta[v^\delta]\|_{L^\infty(\T^d)}^2\big)\dd t.
\end{align*}
Collecting these estimates, we conclude the proof.
\end{proof}

%%%%%%%%%%%%%%%%%%%

\subsection{Passage to the unregularized Hellinger estimate}

We pass to the limit $\eps\to 0$ in the inequality of Lemma \ref{lem.ehi}. 

\begin{lemma}[Unregularized Hellinger estimate]
Let $u$ be an entropy--Fisher weak solution to \eqref{1.eq} and $v^\delta$ be a regularized solution in the sense of Proposition \ref{prop.delta}. Then there exists a constant $C>0$, independent of $\delta$, such that for a.e.\ $\tau\in(0,T]$, 
\begin{align}\label{5.uhi}
  \mathcal H[u|v^\delta](\tau) - \mathcal H[u_0|v_0^\delta]
  \le C\sum_{i=1}^n\int_0^\tau\int_{\T^d}(u_iv_i^\delta)^{1/2}
  |\na(p_i[u]-p_i^\delta[v^\delta])|^2\dx\dd t.
\end{align}
\end{lemma}

\begin{proof}
We drop the nonnegative Hellinger dissipation term on the left-hand side of \eqref{5.ehi} and pass to the limit $\eps\to 0$ in the remaining terms. The limit of the Hellinger distance $\mathcal H_\eps$ follows directly from dominated convergence. The last term in \eqref{5.ehi} vanishes in the limit. It remains to perform the limit in the first term on the right-hand side of \eqref{5.ehi}. Since $(U_i^\eps V_i^\eps)^{1/2}\to (u_iv_i^\delta)^{1/2}$ a.e., the integrand is bounded by an $L^1((0,T)\times\T^d)$ function, taking into account the mass bounds for $u_i$ and $v_i^\delta$ and the $L^2(0,T;L^\infty(\T^d))$ bounds for $\na p_i[u]$ and $\na p_i^\delta[v^\delta]$. This gives
\begin{align*}
  \int_0^\tau\int_{\T^d}
  \sqrt{U_i^\eps V_i^\eps}|\na(p_i[u]-p_i^\delta[v^\delta])|^2
  \dx\dd t \to \int_0^\tau\int_{\T^d}(u_iv_i^\delta)^{1/2}
  |\na(p_i[u]-p_i^\delta[v^\delta])|^2\dx\dd t.
\end{align*}
The proof is finished.
\end{proof}

%%%%%%%%%%%%%%%%%%%%%%%%%%%%%%%%%%%%%%%%%%%%%%%%%%%%%%%%%%%%%%

\section{Interaction energy and modulated stability}\label{sec.ener}

In the previous section, we derived an estimate for the Hellinger distance between an entropy--Fisher weak solution $u$ and a regularized solution $v^\delta$. The main term in estimate \eqref{5.uhi} contains the relative field $\na p_i[u]-\na p_i^\delta[v^\delta]$. Introducing the kernel mismatch
\begin{align}\label{6.Edelta}
  E_i^\delta[u] := \na p_i[u]-\na p_i^\delta[u]
  = \sum_{j=1}^n\na(K_{ij}-K_{ij}^\delta)*u_j,
\end{align} 
the relative field splits into two parts,
\begin{align}\label{6.nap}
  \na(p_i[u]-p_i^\delta[v^\delta])
  = \na(p_i^\delta[u]-p_i^\delta[v^\delta]) 
  + E_i^\delta[u].
\end{align}
The purpose of this section is to control the regularized relative field $\na p_i^\delta[u]-\na p_i^\delta[v^\delta]$ by the regularized interaction energy, while $E_i^\delta[u]$ will be treated as a kernel approximation error. Recall the definitions of the regularized interaction energy and the corresponding dissipation:
\begin{align*}
  \mathcal P^\delta[f] = \frac12\sum_{i=1}^n\int_{\T^d}
  f_i p_i^\delta[f]\dx, \quad
  \mathcal D^\delta[f] = \sum_{i=1}^n\int_{\T^d}
  \na p_i^\delta[f]\cdot\na f_i\dx.
\end{align*}

\subsection{Relative interaction energy identity}

We compute the evolution identity for the regularized relative interaction energy $\mathcal P^\delta[u-v^\delta]$. 

\begin{lemma}[Relative interaction energy identity]
Let $u$ be an entropy--Fisher weak solution to \eqref{1.eq} and $v^\delta$ be a regularized solution in the sense of Proposition \ref{prop.delta}. Then for a.e.\ $\tau\in(0,T]$, 
\begin{align}\label{5.riei}
  \mathcal P^\delta&[u-v^\delta](\tau) 
  - \mathcal P^\delta[u_0-v_0^\delta]
  + \int_0^\tau\mathcal D^\delta[u-v^\delta]\dd t \\
  &\phantom{xx}+ \frac12\sum_{i=1}^n\int_0^\tau\int_{\T^d}
  (u_i+v_i^\delta)|\na(p_i^\delta[u]-p_i^\delta[v^\delta])|^2\dx\dd t 
  \nonumber \\
  &= -\frac12\sum_{i=1}^n\int_0^\tau\int_{\T^d}
  (u_i-v_i^\delta)\na(p_i^\delta[u]+p_i^\delta[v^\delta])\cdot
  \na(p_i^\delta[u]-p_i^\delta[v^\delta])\dx\dd t \nonumber \\
  &\phantom{xx}- \sum_{i=1}^n\int_0^\tau\int_{\T^d}
  u_i E_i^\delta[u]\cdot\na(p_i^\delta[u]-p_i^\delta[v^\delta])
  \dx\dd t. \nonumber 
\end{align}
\end{lemma}

\begin{proof}
Set $w:=u-v^\delta$ and $z:=\rho_\delta*w$. Since $K_{ij}^\delta=\rho_\delta*\rho_\delta*K_{ij}$ and $\rho_\delta$ is even, we have $\mathcal P^\delta[w]=\mathcal P[z]$. We first show that $t\mapsto\mathcal P[z(t)]$ is absolutely continuous. Choosing $\varphi(t,x)=\psi(t)\rho_\delta(y-x)$ with $\psi\in C_c^\infty(0,T)$ and $y\in\T^d$ as a test function in the weak formulations of \eqref{1.eq} and \eqref{4.approx}, we obtain
\begin{align*}
  \pa_t z_i = \na\rho_\delta*\big(F_i[u]-F_i^\delta[v^\delta]\big),
\end{align*}
where the convolution includes the scalar product in $\R^d$. Since $F_i[u]\in L^2(0,T;L^1(\T^d;\R^d))$ by Definition \ref{def.ex} and $F_i^\delta[v^\delta]\in L^2(0,T;L^1(\T^d;\R^d))$ by \eqref{4.est}, Young's convolution inequality gives
\begin{align*}
  \|\pa_t z_i(t)\|_{L^\infty(\T^d)}\le\|\na\rho_\delta\|_{L^\infty(\T^d)}
  \|F_i[u](t)-F_i^\delta[v^\delta](t)\|_{L^1(\T^d)},
\end{align*}
i.e.\ $\pa_t z_i\in L^2(0,T;L^\infty(\T^d))$. In particular, $t\mapsto z_i(t)$ is absolutely continuous with values in $L^1(\T^d)$. Moreover, by mass conservation,
\begin{align*}
  \|z_i(t)\|_{L^\infty(\T^d)}\le\|\rho_\delta\|_{L^\infty(\T^d)}\|w_i(t)\|_{L^1(\T^d)}
  \le\|\rho_\delta\|_{L^\infty(\T^d)}\big(\|u_{i,0}\|_{L^1(\T^d)}+\|v_{i,0}^\delta\|_{L^1(\T^d)}\big)=:M_i^\delta .
\end{align*}
By the symmetry $K_{ij}(x)=K_{ji}(-x)$ (see the proof of Lemma \ref{lem.ec1}),
\begin{align}\label{6.Pdiff}
  \mathcal P[z(t)]-\mathcal P[z(s)]
  = \frac12\sum_{i,j=1}^n\int_{\T^d}\big(z_i(t)+z_i(s)\big)\,
  K_{ij}*\big(z_j(t)-z_j(s)\big)\dx ,
\end{align}
and Young's convolution inequality yields for $0\le s<t\le T$
\begin{align*}
  \big|\mathcal P[z(t)]-\mathcal P[z(s)]\big|
  &\le \sum_{i,j=1}^n M_i^\delta\|K_{ij}\|_{L^1(\T^d)}\|z_j(t)-z_j(s)\|_{L^1(\T^d)} \\
  &\le C_\delta\sum_{j=1}^n\int_s^t\|\pa_t z_j(r)\|_{L^1(\T^d)}\dd r.
\end{align*}
Hence, $t\mapsto\mathcal P^\delta[w(t)]$ is absolutely continuous on $[0,T]$ and, in particular, differentiable a.e.

We identify its derivative. Let $t\in(0,T)$ be a Lebesgue point of $\pa_t z\in L^1(0,T;C(\T^d;\R^n))$. Then $h^{-1}(z(t+h)-z(t))\to\pa_t z(t)$ and $z(t+h)\to z(t)$ in $C(\T^d;\R^n)$ as $h\to 0$. Dividing \eqref{6.Pdiff} (with $t$ and $s$ replaced by $t+h$ and $t$, respectively) by $h$ and passing to the limit $h\to 0$, we infer from Young's convolution inequality and the symmetry of $(K_{ij})$ that
\begin{align*}
  \frac{\dd}{\dd t}\mathcal P^\delta[w(t)]
  = \sum_{i,j=1}^n\int_{\T^d}z_i(t)\,K_{ij}*\pa_t z_j(t)\dx
  = \sum_{i,j=1}^n\int_{\T^d}\pa_t z_i(t)\,K_{ij}*z_j(t)\dx.
\end{align*}
We insert $\pa_t z_i=\na\rho_\delta*G_i$ and use that $\na\rho_\delta$ is odd, since $\rho_\delta$ is even. Together with the representation of $p_i^\delta[w]$ from above, this gives
\begin{align*}
  \frac{\dd}{\dd t}\mathcal P^\delta[w]
  = -\sum_{i=1}^n\int_{\T^d}G_i\cdot\na\bigg(\rho_\delta*\sum_{j=1}^n
  K_{ij}*z_j\bigg)\dx
  = -\sum_{i=1}^n\int_{\T^d}G_i\cdot\na p_i^\delta[w]\dx.
\end{align*}
Since $G_i=\na w_i+u_i\na p_i[u]-v_i^\delta\na p_i^\delta[v^\delta]$, we conclude that
\begin{align*}
  \frac{\dd}{\dd t}\mathcal P^\delta[w] 
  = -\sum_{i=1}^n\int_{\T^d}\na p_i^\delta[w]\cdot\na w_i\dx 
  - \sum_{i=1}^n\int_{\T^d}\na p_i^\delta[w]\cdot
  \big(u_i\na p_i[u] - v_i^\delta\na p_i^\delta[v^\delta]\big)\dx.
\end{align*}
The first term equals $-\mathcal D^\delta[w]$. For the drift part, we use $\na p_i[u] = \na p_i^\delta[u] + E_i^\delta[u]$, recalling definition \eqref{6.Edelta} of $E_i^\delta[u]$. We infer that
\begin{align*}
  u_i&\na p_i[u] - v_i^\delta\na p_i^\delta[v^\delta]
  = u_i\na p_i^\delta[u] - v_i^\delta\na p_i^\delta[v^\delta]
  + u_iE_i^\delta[u] \\
  &= \frac12(u_i+v_i^\delta)\na(p_i^\delta[u]-p_i^\delta[v^\delta])
  + \frac12(u_i-v_i^\delta)\na(p_i^\delta[u]+p_i^\delta[v^\delta])
  + u_iE_i^\delta[u].
\end{align*}
We substitute this identity and integrate over $(0,\tau)$ to conclude \eqref{5.riei}.
\end{proof}

%%%%%%%%%%%%%%%%%

\subsection{Estimated interaction inequality}

We estimate the relative interaction energy identity \eqref{5.riei}. We abbreviate
\begin{align*}
  \Sigma(t) := \sum_{i=1}^n
  \big(\|\na p_i^\delta[u]\|_{L^\infty(\T^d)}^2
  + \|\na p_i^\delta[v^\delta]\|_{L^\infty(\T^d)}^2\big).
\end{align*}

\begin{lemma}[Estimated interaction inequality]
Let $u$ be an entropy--Fisher weak solution to \eqref{1.eq} and $v^\delta$ be a regularized solution in the sense of Proposition \ref{prop.delta}. Then there exists a constant $C>0$, independent of $\delta$, such that for a.e.\ $\tau\in(0,T]$,
\begin{align}\label{6.eii}
  \mathcal P^\delta&[u-v^\delta](\tau)
  - \mathcal P^\delta[u_0-v_0^\delta]
  + \frac12\int_0^\tau\mathcal D^\delta[u-v^\delta]\dd t \\
  &\phantom{xx}+ \frac14\sum_{i=1}^n\int_0^\tau\int_{\T^d}
  (u_i+v_i^\delta)
  |\na(p_i^\delta[u]-p_i^\delta[v^\delta])|^2\dx\dd t \nonumber \\
  &\le C\int_0^\tau\Sigma(t)\mathcal H[u|v^\delta](t)\dd t
  + C\sum_{i=1}^n\int_0^\tau\|E_i^\delta[u(t)]\|_{L^\infty(\T^d)}^2
  \dd t. \nonumber 
\end{align}
\end{lemma}

\begin{proof}
We start from the relative interaction energy identity \eqref{5.riei}. We estimate the first term on the right-hand side of \eqref{5.riei}. We infer from $u_i-v_i^\delta = (u_i^{1/2}-(v_i^\delta)^{1/2})(u_i^{1/2}+(v_i^\delta)^{1/2})$ and the Young inequality that
\begin{align*}
  \bigg|\frac12&\sum_{i=1}^n\int_0^\tau\int_{\T^d}
  (u_i-v_i^\delta)\na(p_i^\delta[u]+p_i^\delta[v^\delta])\cdot
  \na(p_i^\delta[u]-p_i^\delta[v^\delta])\dx\dd t\bigg| \\
  &\le \frac{1}{16}\sum_{i=1}^n\int_0^\tau\int_{\T^d}
  \big(u_i^{1/2}+(v_i^\delta)^{1/2}\big)^2
  \big|\na(p_i^\delta[u]-p_i^\delta[v^\delta])\big|^2\dx\dd t \\
  &\phantom{xx}+ \sum_{i=1}^n\int_0^\tau\int_{\T^d}
  \big(u_i^{1/2}-(v_i^\delta)^{1/2}\big)^2
  \big|\na(p_i^\delta[u]+p_i^\delta[v^\delta])\big|^2\dx\dd t \\
  &\le \frac18\sum_{i=1}^n\int_0^\tau\int_{\T^d}(u_i+v_i^\delta)
  \big|\na(p_i^\delta[u]-p_i^\delta[v^\delta])\big|^2\dx\dd t
  + C\int_0^\tau \Sigma(t)\mathcal H[u|v^\delta](t)\dd t.
\end{align*}
It remains to estimate the kernel mismatch term. Again by the Young inequality,
\begin{align*}
  \bigg|\sum_{i=1}^n&\int_0^\tau\int_{\T^d}
  u_i E_i^\delta[u]\cdot\na(p_i^\delta[u]-p_i^\delta[v^\delta])
  \dx\dd t\bigg| \\
  &\le \frac18\sum_{i=1}^n\int_0^\tau\int_{\T^d}u_i
  \big|\na(p_i^\delta[u]-p_i^\delta[v^\delta])\big|^2\dx\dd t
  + 2\sum_{i=1}^n\int_0^\tau\int_{\T^d}u_i
  |E_i^\delta[u]|^2\dx\dd t.
\end{align*}
Substituting the previous estimates into \eqref{5.riei} and taking into account the mass bound for $u_i$ concludes the proof.
\end{proof}

%%%%%%%%%%%%%%%%

\subsection{Regularized weak--strong stability}

We combine the Hellinger estimate with the interaction inequality. Define, for some $\lambda>0$ to be defined later,
\begin{align*}
  \mathcal M^\delta[u|v^\delta] = \mathcal H[u|v^\delta]
  + \lambda\mathcal P^\delta[u-v^\delta].
\end{align*}

\begin{proposition}[Regularized weak--strong stability]
\label{prop.rwss}
Let $u$ be an entropy--Fisher weak solution to \eqref{1.eq} and $v^\delta$ be a regularized solution in the sense of Proposition \ref{prop.delta}. Then there exist constants $c>0$ and $C_T>0$ depending on $T$, independent of $\delta$, such that for a.e.\ $\tau\in(0,T]$,
\begin{align}\label{6.rwss}
  \mathcal M^\delta&[u|v^\delta](\tau)
  + c\int_0^\tau\mathcal D^\delta[u-v^\delta]\dd t \\
  &\phantom{xx}+ c\sum_{i=1}^n\int_0^\tau\int_{\T^d}
  (u_i+v_i^\delta)\big|\na(p_i^\delta[u]-p_i^\delta[v^\delta])\big|^2
  \dx\dd t \nonumber \\
  &\le C_T\bigg(\mathcal M^\delta[u_0|v_0^\delta]
  + \sum_{i=1}^n\int_0^\tau\|E_i^\delta[u]\|_{L^\infty(\T^d)}^2
  \dd t\bigg). \nonumber 
\end{align}
\end{proposition}

\begin{proof}
We estimate the right-hand side of \eqref{5.uhi}. It follows from the decomposition \eqref{6.nap} and Young's inequality that
\begin{align*}
  (u_iv_i^\delta)^{1/2}\big|\na(p_i[u]-p_i^\delta[v^\delta])\big|^2
  \le (u_i+v_i^\delta)\big|\na(p_i^\delta[u]-p_i^\delta[v^\delta])\big|^2
  + (u_i+v_i^\delta)|E_i^\delta[u]|^2.
\end{align*}
Since both $u_i$ and $v_i^\delta$ conserve mass, \eqref{5.uhi} becomes
\begin{align}\label{6.aux}
  \mathcal H&[u|v^\delta](\tau) - \mathcal H[u_0|v_0^\delta] \\
  &\le C\sum_{i=1}^n\int_0^\tau\int_{\T^d}(u_i+v_i^\delta)
  \big|\na(p_i^\delta[u]-p_i^\delta[v^\delta])\big|^2\dx\dd t
  +  C\sum_{i=1}^n\int_0^\tau\|E_i^\delta[u]\|_{L^\infty(\T^d)}^2\dd t.
  \nonumber 
\end{align}
Next, we multiply \eqref{6.eii} by $\lambda>0$ and add the resulting inequality to \eqref{6.aux}. Choosing $\lambda$ sufficiently large, the first term on the right-hand side of \eqref{6.aux} is absorbed by the last term of the left-hand side of $\lambda\times$\eqref{6.eii}. Thus, for some $c>0$,
\begin{align}\label{6.aux2} 
  \mathcal M^\delta&[u|v^\delta](\tau) 
  + c\int_0^\tau\mathcal D^\delta[u-v^\delta]\dd t \\
  &\phantom{xx}+ c\sum_{i=1}^n\int_0^\tau\int_{\T^d}(u_i+v_i^\delta)
  \big|\na(p_i^\delta[u]-p_i^\delta[v^\delta])\big|^2\dx\dd t 
  \nonumber \\
  &\le \mathcal M^\delta[u_0|v_0^\delta]
  + C\int_0^\tau\Sigma\mathcal H[u|v^\delta]\dd t
  + C\sum_{i=1}^n\int_{0}^\tau\|E_i^\delta[u]\|_{L^\infty(\T^d)}^2\dd t.
  \nonumber 
\end{align}
Since the regularized matrix kernel is positive semidefinite, we have $\mathcal H[u|v^\delta]\le\mathcal M^\delta[u|v^\delta]$, and by the endpoint drift estimate in Lemma \ref{lem.ec2}, it holds that $\int_0^T\Sigma\dd t\le C_T$. Therefore, we can apply Gronwall's lemma to \eqref{6.aux2} to obtain
\begin{align*}
  \mathcal M^\delta[u|v^\delta](\tau)
  \le C_T\bigg(\mathcal M^\delta[u_0|v_0^\delta]
  + \sum_{i=1}^n\int_0^\tau\|E_i^\delta[u]\|_{L^\infty(\T^d)}^2
  \dd t\bigg).
\end{align*}
Substituting this bound back into \eqref{6.aux2} yields \eqref{6.rwss}.
\end{proof}

%%%%%%%%%%%%%%%%%%%

\subsection{Proof of Theorem \ref{thm.stab}}

Now we can pass to the limit $\delta\to 0$ and prove the weak--weak stability and uniqueness. For this, let $u$ be an entropy--Fisher weak solution with initial datum $u_0$, and let $v_0$ be another admissible initial datum. Choose smooth strictly positive initial data $v_0^\delta$ with the same componentwise masses as $v_0$, such that $v_0^\delta\to v_0$ in $L^1(\T^d;\R^n)$ and with uniform $L\log L(\T^d;\R^n)$ bounds. Then $\mathcal M^\delta[u_0|v_0^\delta]\to\mathcal M[u_0|v_0]$. Indeed, the Hellinger part is continuous under $L^1(\T^d)$ convergence because $(\sqrt{a}-\sqrt{b})^2\le|a-b|$. We deduce from $\mathcal P^\delta[u_0-v_0^\delta] = \mathcal P[\rho_\delta*(u_0-v_0^\delta)]$ and the continuity of $\mathcal P[\cdot]$ under $L^1(\T^d)$ convergence by Lemma \ref{lem.ec1} that $\mathcal P^\delta[u_0-v_0^\delta]\to\mathcal P[u_0-v_0]$. 

Let $v^\delta$ be a regularized solution to \eqref{4.approx} with initial datum $v_0^\delta$. Proposition \ref{prop.rwss} gives
\begin{align}\label{6.rwss2}
  \mathcal M^\delta[u|v^\delta](\tau)
  + c\int_0^\tau\mathcal D^\delta[u-v^\delta]\dd t
  \le C_T\bigg(\mathcal M^\delta[u_0|v_0^\delta]
  + \sum_{i=1}^n\int_0^\tau\|E_i^\delta[u]\|_{L^\infty(\T^d)}^2
  \dd t\bigg).  
\end{align}
We show that the last term vanishes in the limit $\delta\to 0$. It follows from the definition $K_{ij}^\delta = \eta_\delta*K_{ij}$ and the endpoint drift estimate in Lemma \ref{lem.ec2}, applied to $u_j-\eta_\delta*u_j$, that 
\begin{align*}
  \|E_i^\delta[u]\|_{L^\infty(\T^d)} 
  &= \|\na(p_i[u] - p_i^\delta[u])\|_{L^\infty(\T^d)}
  = \bigg\|\sum_{j=1}^n \na K_{ij} * (u_j-\eta_\delta * u_j)
  \bigg\|_{L^\infty(\T^d)} \\
  &\le C\sum_{j=1}^n\|\na K_{ij}\|_{L^{d,\infty}(\T^d)}
  \|\na u_j - \eta_\delta*\na u_j\|_{L^1(\T^d)}. 
\end{align*}
Thus, since $\eta_\delta*\na u_j\to \na u_j$ strongly in $L^2(0,T;L^1(\T^d))$,
\begin{align*}
  \sum_{i=1}^n\int_0^T\|E_i^\delta[u]\|_{L^\infty(\T^d)}^2\dd t
  \to 0.
\end{align*}

We pass to the limit $\delta\to 0$ in \eqref{6.rwss2}. We know from Theorem \ref{thm.ex} that, up to a subsequence, $v^\delta(t)\to \bar{v}(t)$ strongly in $L^1(\T^d)$ for a.e.\ $t\in(0,T)$. The limit $\bar v$ is an entropy--Fisher weak solution to \eqref{1.eq} with initial datum $v_0$. Arguing similarly as for the initial data, we have $\mathcal M^\delta[u|v^\delta](t)\to\mathcal M[u|\bar v](t)$. We pass to the limit in the dissipation term by exploiting the lower semicontinuity of the Fourier quadratic term. For this, we set $f^\delta=u-v^\delta$ and $f=u-\bar v$. The strong $L^1(\T^d)$ convergence gives $\widehat{f}^\delta(t,\xi)\to\widehat{f}(t,\xi)$ for every $\xi\in\mathbb Z^d$ and a.e.\ $t\in[0,T]$. Moreover, for every $\xi\in\mathbb Z^d$, $\widehat{K}^\delta(\xi)=|\widehat{\rho}_\delta(\xi)|^2 \widehat{K}(\xi)$, $\widehat{\rho}_\delta(\xi) = \mathcal F_{\R^d}[\rho](\delta\xi)\to 1$ as $\delta\to 0$, where $\rho_\delta$ is the periodization of $\delta^{-d}\rho(\cdot/\delta)$, and $\widehat{K}(\xi)$ is positive semidefinite. Since
\begin{align*}
  \mathcal D^\delta[f^\delta](t) = 4\pi^2\sum_{\xi\in\mathbb Z^d}|\xi|^2\widehat{f^\delta}(t,\xi)^* \widehat{K^\delta}(\xi) \widehat{f^\delta}(t,\xi)\ge0,
\end{align*}
where $\widehat{f^\delta}(t,\xi)^*$ denotes the conjugate transpose, Fatou's lemma in $(t,\xi)$ leads to
\begin{align*}
  \int_0^{\tau}\mathcal D[u-\bar v]\dd t
  \le \liminf_{\delta\to 0}\int_0^{\tau}
  \mathcal D^\delta[u-v^\delta]\dd t.
\end{align*}
Combining these results, we obtain for a.e.\ $\tau\in(0,T)$,
\begin{align}\label{6.aux3}
  \mathcal M[u|\bar v](\tau) + c\int_0^{\tau}\mathcal D[u-\bar v]\dd t
  \le C_T\mathcal M[u_0|v_0]. 
\end{align}
If $u_0=v_0$, the right-hand side vanishes, implying that $\mathcal M[u|\bar v](\tau)=0$ for a.e.\ $\tau\in(0,T)$ and therefore $u=\bar v$ a.e. This proves that entropy--Fisher weak solutions are unique for each admissible initial datum. Finally, let $v$ be any entropy--Fisher weak solution with initial datum $v_0$. By the uniqueness just proved, $v=\bar v$, and we obtain \eqref{6.aux3} for $v$ instead of $\bar v$. 

%%%%%%%%%%%%%%%%%%%%%%%%%%%%%%%%%%%%%%%%%%%%%%%%%%%%%%%%%%%%%%

\section{Positive-time $L^2(\T^d)$ regularization and convergence to equilibrium}\label{sec.regul}

In this section, we prove Theorem \ref{thm.time}. We start with a regularized solution $u^\delta$ and then pass to the limit $\delta\to 0$.

\subsection{Weighted $L^2(\T^d)$ regularization and decay for the regularized solutions}

The endpoint framework developed above yields more than well-posedness.
We derive a quantitative $L^2(\T^d)$ regularization and decay estimate for the regularized solutions. 

\begin{lemma}[Exponential decay of the regularized relative entropy]
\label{lem.regexp}
Let $u^\delta$ be the regularized solution to \eqref{4.approx} from Proposition \ref{prop.delta}. Then
\begin{align*}
  \mathcal E[u^\delta(t)|m] \le e^{-\kappa t}\mathcal E[u_0^\delta|m]
  \quad\mbox{for }t\ge 0,
\end{align*}
where $\kappa= 4/C_{LS}$, $C_{LS}>0$ is the constant from the logarithmic Sobolev inequality in Lemma \ref{lem.lsi}, and $\mathcal E$ is the relative entropy defined in \eqref{4.E}. Moreover, for $t\ge 0$,
\begin{align*}%\label{7.L1}
  \sum_{i=1}^n\frac{1}{m_i}\|u_i^\delta(t)-m_i\|_{L^1(\T^d)}^2
  \le 2e^{-\kappa t}\mathcal{E}[u_0^\delta|m].
\end{align*}
\end{lemma}

\begin{proof}
We sum the logarithmic Sobolev inequality
\begin{align*}
  \int_{\T^d}u_i^\delta\log\frac{u_i^\delta}{m_i}\dx
  \le C_{LS}\int_{\T^d}|\na(u_i^\delta)^{1/2}|^2\dx
\end{align*}
over $i=1,\ldots,n$ and substitute the resulting expression into the entropy identity \eqref{4.dEdt}, which gives 
\begin{align*}
  \frac{\dd}{\dd t}\mathcal E[u^\delta|m] 
  + \kappa\mathcal E[u^\delta|m]\le 0.
\end{align*}
Gronwall's lemma shows the entropy decay, and the Csisz\'ar--Kullback--Pinsker inequality in Lemma \ref{lem.ckpi}, applied to the probability density $u_i^\delta/m_i$, implies that
\begin{align*}
  \|u_i^\delta-m_i\|_{L^1(\T^d)}^2
  \le 2m_i\int_{\T^d}u_i^\delta\log\frac{u_i^\delta}{m_i}\dx.
\end{align*}
The claimed estimate follows after dividing this inequality by $m_i$ and summing over $i=1,\ldots,n$. 
\end{proof}

Let $w_i^\delta=u_i^\delta-m_i$ and $w^\delta=(w_1^\delta,\ldots,w_n^\delta)$. For $p=1,2$, define the weighted squared $L^p(\T^d)$ norms
$$
  \mathscr L_1^\delta(t) := \sum_{i=1}^n\frac1{m_i}
  \|w_i^\delta(t)\|_{L^1(\T^d)}^2, \quad 
  \mathscr L_2^\delta(t) := \sum_{i=1}^n\frac1{m_i}
  \|w_i^\delta(t)\|_{L^2(\T^d)}^2.
$$

\begin{lemma}[Weighted $L^2(\T^d)$ regularization and decay for $u^\delta$]\label{lem.L2}
Let $\mathcal E_0^\delta:=\mathcal E[u_0^\delta|m]$. Then
\begin{align*}
  \mathscr L_2^\delta(t) \le 2\mathcal E_0^\delta
  (\kappa C_{\rm N})^{d/2}\exp(C_K\mathcal E_0^\delta)
  (e^{2\kappa t/d}-1)^{-d/2}, \quad t>0,
\end{align*}
where $C_K = C_d\sum_{i,j=1}^n m_j\|\na K_{ij}\|_{L^{d,\infty}(\T^d)}^2$ and $C_N>0$ is the constant of the Nash inequality in Lemma \ref{lem.nash}.
\end{lemma}

\begin{proof}
Since the kernels have zero mean, the constant vector $m$ does not contribute to the drift field, i.e.\ $\na p_i^\delta[u^\delta] = \na p_i^\delta[w^\delta]$. Hence,
\begin{align*}
  \pa_t w_i^\delta = \Delta w_i^\delta 
  + \diver\big((w_i^\delta+m_i)\na p_i^\delta[w^\delta]\big).
\end{align*}
We multiply this equation by $w_i^\delta/m_i$, sum over $i=1,\ldots,n$, and use Young's inequality:
\begin{align*}
  & \frac12\frac{\dd}{\dd t}\mathscr L_2^\delta(t)
  + \sum_{i=1}^n\frac1{m_i}\|\na w_i^\delta\|_{L^2(\T^d)}^2
  + \mathcal D^\delta[w^\delta]
  \le \frac12\sum_{i=1}^n\frac1{m_i}\|\na w_i^\delta\|_{L^2(\T^d)}^2
  + \frac12 B_\delta(t)\mathscr L_2^\delta(t), \\
  &\mbox{where }B_\delta(t) = \sum_{i=1}^n
  \|\na p_i^\delta[u^\delta](t)\|_{L^\infty(\T^d)}^2.
\end{align*}
We absorb the gradient term on the right-hand side and drop the term $\mathcal D^\delta[w^\delta]$:
\begin{align}\label{7.dLdt}
  \frac{\dd}{\dd t}\mathscr L_2^\delta(t)
  + \sum_{i=1}^n\frac1{m_i}\|\na w_i^\delta(t)\|_{L^2(\T^d)}^2
  \le B_\delta(t)\mathscr L_2^\delta(t).
\end{align}
By the zero-mean Nash inequality in Lemma \ref{lem.nash}, since $\int_{\T^d}w_i^\delta\dx=0$,
\begin{align*}
  \frac{1}{m_i}\|\na w_i^\delta\|_{L^2(\T^d)}^2
  \ge \frac{1}{C_N}\frac{(m_i^{-1}\|w_i^\delta
  \|_{L^2(\T^d)}^2)^{1+2/d}}{(m_i^{-1}\|w_i^\delta
  \|_{L^1(\T^d)}^2)^{2/d}}.
\end{align*}
We sum this inequality over $i=1,\ldots,n$ and use the H\"older inequality with the conjugate exponents $p=(1+c)/c$ and $p'=1+c$, where $c=2/d$, in the form
\begin{align*}
  \sum_{i=1}^n\frac{a_i^{1+2/d}}{b_i^{2/d}}
  \ge \frac{(\sum_{i=1}^n a_i)^{1+2/d}}{(\sum_{i=1}^n b_i)^{2/d}}
  \quad\mbox{for }a_i,\,b_i>0,
\end{align*}
to infer that
\begin{align*}
  \sum_{i=1}^n\frac1{m_i}\|\na w_i^\delta(t)\|_{L^2(\T^d)}^2
  \ge \frac{1}{C_N}\frac{(\mathscr L_2^\delta(t))^{1+2/d}}{
  (\mathscr L_1^\delta(t))^{2/d}}.
\end{align*}
It follows from Lemma \ref{lem.regexp} that $\mathscr L_1^\delta(t) \le 2e^{-\kappa t}\mathcal E_0^\delta$. Then \eqref{7.dLdt} leads to
\begin{align*}
  \frac{\dd}{\dd t}\mathscr L_2^\delta(t)
  + C_N^{-1}(2\mathcal E_0^\delta)^{-2/d} e^{2\kappa t/d}
  (\mathscr L_2^\delta(t))^{1+2/d} \le B_\delta(t)
  \mathscr L_2^\delta(t).
\end{align*}
We solve this differential inequality. To this end, let $\Gamma_\delta(t) = \int_0^t B_\delta(s)\dd s\ge 0$. We multiply the inequality  by $e^{-\Gamma_\delta(t)}$:
\begin{align*}
  \frac{\dd}{\dd t}\big(e^{-\Gamma_\delta(t)}\mathscr L_2^\delta\big)
  + C_N^{-1}(2\mathcal E_0^\delta)^{-2/d} e^{2\kappa t/d}
  e^{-\Gamma_\delta(t)}(\mathscr L_2^\delta\big)^{1+2/d} \le 0,
\end{align*}
giving
\begin{align*}
  \frac{\dd}{\dd t}\big(e^{-\Gamma_\delta(t)}\mathscr L_2^\delta
  \big)^{-2/d} \ge \frac{2}{d}\frac{1}{C_N}
  (2\mathcal E_0^\delta)^{-2/d} e^{2\kappa t/d}e^{2\Gamma_\delta/d}.
\end{align*}
An integration in time yields, after a computation, that
\begin{align*}
  \mathscr L_2^\delta(t) \le 2\mathcal E_0^\delta
  (\kappa C_N)^{d/2}\exp\bigg(\int_0^t B_\delta(s)\dd s\bigg)
  \big(e^{2\kappa t/d}-1\big)^{-d/2}.
\end{align*}
We estimate the drift integral $\int_0^t B_\delta(s)\dd s$. We deduce from the endpoint drift estimate in Lemma \ref{lem.ec2} and the first estimate in Corollary \ref{coro}, applied to the regularized drift $\na p_i^\delta[u^\delta] = \na p_i[\eta_\delta* u^\delta]$, and from the uniform estimate from Proposition \ref{prop.delta} that
\begin{align*}
  \int_0^\infty B_\delta(s)\dd s\le C_d\sum_{i,j=1}^n m_j
  \|\na K_{ij}\|_{L^{d,\infty}(\T^d)}^2\mathcal E_0^\delta.
\end{align*}
The lemma follows after substitution into the previous bound.
\end{proof}

%%%%%%%%%%%%%%%%

\subsection{Proof of Theorem \ref{thm.time}}

We establish the positive-time regularization and decay estimates for the entropy--Fisher weak solution by passing to the limit from the corresponding estimates for the regularized solutions. This
gives the existence of an entropy--Fisher weak solution with the same regularization and convergence-to-equilibrium properties as the approximate solutions. Let $u_0^\delta$ be the regular initial approximation satisfying \eqref{4.init}, and let $u^\delta$ be the corresponding regularized solution from Proposition \ref{prop.delta}. By Lemma \ref{lem.L2} and the uniform bounds on the initial data \eqref{4.init}, there exists a constant $C>0$ independent of $\delta$, such that
\begin{align*}
  \sum_{i=1}^n\frac{1}{m_i}\|u_i^\delta(t)-m_i\|_{L^2(\T^d)}^2
  \le C\big(e^{2\kappa t/d}-1\big)^{-d/2}. 
\end{align*}
By the compactness argument in Theorem~\ref{thm.ex}, we have, up to a subsequence, $u_i^\delta(t)\to u_i(t)$ strongly in $L^1(\T^d)$ for a.e.\ $t>0$, where $u$ is an entropy--Fisher weak solution with initial datum $u_0$. The uniform $L^2(\T^d)$ bound above implies that $u_i^\delta(t)-m_i$ is bounded in $L^2(\T^d)$. Hence,
\begin{align*}
  u_i^\delta(t)-m_i \rightharpoonup u_i(t)-m_i 
  \quad\text{weakly in }L^2(\T^d) 
\end{align*} 
for a.e.\ $t>0$. By weak lower semicontinuity,
\begin{align*}%\label{7.aux} 
  \sum_{i=1}^n\frac{1}{m_i}\|u_i(t)-m_i\|_{L^2(\T^d)}^2
  \le \liminf_{\delta\to 0}
  \sum_{i=1}^n\frac{1}{m_i}\|u_i^\delta(t)-m_i\|_{L^2(\T^d)}^2.
\end{align*}
This proves that the limit entropy--Fisher weak solution $u$ satisfies the regularization estimate 
\begin{align}\label{7.reg}
  \sum_{i=1}^n\frac1{m_i}\|u_i(t)-m_i\|_{L^2(\T^d)}^2
  \le C\big(e^{2\kappa t/d}-1\big)^{-d/2}.
\end{align}
By the uniqueness of entropy--Fisher weak solutions (see Theorem \ref{thm.stab}), the constructed limit $u$ coincides with any entropy--Fisher weak solution with initial datum $u_0$; hence \eqref{7.reg} holds for all entropy--Fisher weak solutions.

%%%%%%%%%%%%%%%%%%%%%%%%%%%%%%%%%%%%%%%%%%%%%%%%%

\section{Bounded interaction kernels}
\label{sec:bounded-kernels}

We prove Theorem~\ref{thm:bounded-kernels} by establishing
the elementary estimates for the drift field and interaction energy under the bounded-kernel assumption.
Once they are available, the arguments for existence, weak--weak stability, and $L^2$ regularization as well as decay carry over with only minor modifications.

\subsection{Proof of Theorem~\ref{thm:bounded-kernels}}

The following elementary estimates replace the endpoint bounds of
Lemmas~\ref{lem.ec1} and
\ref{lem.ec2}.

\begin{lemma}[Bounds for the force and interaction energy]
%\label{lem:bounded-kernel-estimates}
Under the kernel assumptions of Theorem~\ref{thm:bounded-kernels},
for $f\in W^{1,1}(\T^d;\R^n)$ and $i=1,\ldots,n$,
\begin{equation}\label{eq:bounded-force}
\|\nabla p_i[f]\|_{L^\infty(\T^d)}
\le \sum_{j=1}^n\|K_{ij}\|_{L^\infty(\T^d)}
\|\nabla f_j\|_{L^1(\T^d)}.
\end{equation}
Furthermore, for all $f,g\in L^1(\T^d;\R^n)$,
\begin{align}
|\mathcal P[f]|
&\le\frac12\sum_{i,j=1}^n
\|K_{ij}\|_{L^\infty(\T^d)}\|f_i\|_{L^1(\T^d)}\|f_j\|_{L^1(\T^d)},
\label{eq:bounded-energy}\\
|\mathcal P[f]-\mathcal P[g]|
&\le\frac12\sum_{i,j=1}^n
\|K_{ij}\|_{L^\infty(\T^d)}
\|f_i+g_i\|_{L^1(\T^d)}\|f_j-g_j\|_{L^1(\T^d)}.
\label{eq:bounded-energy-continuity}
\end{align}
In particular, the interaction energy $\mathcal P$, defined in \eqref{eq:def-interaction-energy}, is locally Lipschitz continuous and nonnegative on $L^1(\T^d;\R^n)$. 
\end{lemma}

\begin{proof}
Since $\nabla(K_{ij}*f_j)=K_{ij}*\nabla f_j$, inequality
\eqref{eq:bounded-force} follows from Young's convolution inequality.
The same inequality yields
$\|K_{ij}*f_j\|_{L^\infty(\T^d)}\le
\|K_{ij}\|_{L^\infty(\T^d)}\|f_j\|_{L^1(\T^d)}$,
which proves \eqref{eq:bounded-energy}. The symmetry
$K_{ij}(x)=K_{ji}(-x)$ gives
\[
\mathcal P[f]-\mathcal P[g]
=\frac12\sum_{i,j=1}^n\int_{\T^d}
(f_i+g_i)\big(K_{ij}*(f_j-g_j)\big)\,\dd x,
\]
and hence \eqref{eq:bounded-energy-continuity}. The nonnegativity of $\mathcal{P}$ follows for $L^2$ functions. In the general case, we approximate $f\in L^1(\T^d)$ by $f_k\in L^2(\T^d)$  with $\mathcal{P}[f_k]\ge 0$, i.e.\ $f_k\to f$ in $L^1(\T^d)$. Then \eqref{eq:bounded-energy-continuity} shows that $\mathcal{P}[f_k]\to\mathcal{P}[f]$ and consequently $\mathcal{P}[f]\ge 0$.
\end{proof}

\begin{proof}[Proof of Theorem~\ref{thm:bounded-kernels}]
By \eqref{eq:bounded-force}--\eqref{eq:bounded-energy-continuity},
the force and flux estimates of
Corollary~\ref{coro}
remain valid.
For the existence proof in Section~\ref{sec.ex},
we use the same regularization
$K_{ij}^\delta=\eta_\delta*K_{ij}$.
The entropy estimates and compactness argument still hold. The Hellinger estimates in
Section~\ref{sec.hell}
apply directly. In
Section~\ref{sec.ener},
the continuity estimate
\eqref{eq:bounded-energy-continuity}
replaces the $L\log L$--BMO argument.
In particular, it ensures the absolute continuity
of $t\mapsto\mathcal P^\delta[u(t)-v^\delta(t)]$,
since
$\mathcal P^\delta[w]=\mathcal P[\rho_\delta*w]$
and $t\mapsto\rho_\delta*w(t)$
is absolutely continuous in $L^1(\T^d)$.
Furthermore, the approximation error \eqref{6.Edelta} satisfies
\begin{equation*}%\label{eq:bounded-kernel-mismatch}
\begin{aligned}
\|E_i^\delta[u]\|_{L^2(0,T;L^\infty(\T^d))}
&\le\sum_{j=1}^n\|K_{ij}\|_{L^\infty(\T^d)}
\|\nabla u_j-\eta_\delta*\nabla u_j\|_{L^2(0,T;L^1(\T^d))}
\longrightarrow0.
\end{aligned}
\end{equation*}
Hence, the interaction energy identity
and the stability argument remain valid,
proving the assertion of
Theorem~\ref{thm.stab}. Finally, the proof of the positive-time $L^2$
regularization and decay in
Section~\ref{sec.regul}
applies without further changes except that the constant $C_K$ in Lemma \ref{lem.L2} has to be replaced by $C\sum_{i,j}m_j\|K_{ij}\|_{L^\infty(\T^d)}^2$.
\end{proof}

\subsection{Examples}

\emph{Kernels of convolution type.}
Let $w\in L^2(\T^d)$, set $\widetilde w(x):=w(-x)$, and define
\begin{align*}
  k := w*\widetilde w - \bigg(\int_{\T^d}w\dx\bigg)^2.
\end{align*}
By the Cauchy--Schwarz inequality, $\|w*\widetilde w\|_{L^\infty(\T^d)}\le\|w\|_{L^2(\T^d)}^2$, so $k\in L^\infty(\T^d)$. Moreover, $k$ is even and has zero mean, and its Fourier coefficients are given by $\widehat k(0)=0$ and $\widehat k(\xi)=|\widehat w(\xi)|^2\ge 0$ for $\xi\in\mathbb Z^d\setminus\{0\}$. Let $A=(a_{ij})_{i,j=1}^n$ be a real symmetric positive semidefinite matrix and set $K_{ij}:=a_{ij}k$. Then $K_{ij}(x)=K_{ji}(-x)$, $\int_{\T^d}K_{ij}\dx=0$, and Parseval's identity yields for all $f\in L^2(\T^d;\R^n)$,
\begin{align*}
  \sum_{i,j=1}^n\int_{\T^d}f_i(K_{ij}*f_j)\dx
  = \sum_{\xi\in\mathbb Z^d}\widehat k(\xi)\,
  \widehat f(\xi)^*A\widehat f(\xi)\ge 0,
\end{align*}
recalling that $\widehat f(\xi)^*$ is the conjugate transpose of $\widehat f(\xi)$. Thus, the kernels $K_{ij}$ satisfy the assumptions of Theorem~\ref{thm:bounded-kernels}. We note that subtracting the mean value does not change the drift field, since $\na p_i[u]$ depends only on $\na K_{ij}$.

For example, let $d=1$, $\T=\R/\mathbb Z$, and $w=\mathbf 1_{(-R,R)}$ with $R\in(0,1/4]$. The kernel $\mathbf 1_{(-R,R)}-2R$ itself is not positive semidefinite, since its Fourier coefficients $\sin(2\pi R\xi)/(\pi\xi)$, $\xi\ne 0$, change sign. However, $\mathbf 1_{(-R,R)}*\mathbf 1_{(-R,R)}$ is the periodic extension of the tent function $(2R-|x|)_+$ for $|x|\le 1/2$, which has the mean value $4R^2$. Therefore, the kernel
\begin{align*}
  k(x) = (2R-|x|)_+ - 4R^2, \quad |x|\le\frac12,
\end{align*}
(extended periodically) is bounded, even, has zero mean, and is positive semidefinite.

\emph{One-dimensional Coulomb kernel.}
On $\T=\R/\mathbb Z$, consider the zero-mean periodic Green function
\[
G(x)=\frac12\left(x^2-x+\frac16\right),\qquad 0\le x<1,
\]
which satisfies $-G''=\delta_0-1$ in the sense of distributions.
Since $G\in L^\infty(\T)$, the matrix kernel
$K_{ij}=z_i z_j G$, with $z_i\in\R$, satisfies the boundedness and
symmetry conditions. It is positive semidefinite
because, writing $\rho[f]=\sum_i z_i f_i$ and
$\phi[f]=G*\rho[f]$,
\[
\sum_{i,j=1}^n\int_{\T}f_i(K_{ij}*f_j)\,\dd x
=\int_{\T}|\partial_x\phi[f]|^2\,\dd x\ge0.
\]
Thus, Theorem~\ref{thm:bounded-kernels} applies to the
one-dimensional Poisson--Nernst--Planck system.

%%%%%%%%%%%%%%%%%%%%%%%%%%%%%%%%%%%%%%%%%%%%%%%%%%%%%%%%%%%%%%

\end{document}